\documentclass{article}

\usepackage{xcolor}
\usepackage{booktabs}
\usepackage[margin=.5cm,labelfont=bf]{caption} 
\usepackage[top=1.8cm, bottom=2.8cm]{geometry}

\usepackage{enumerate}

\usepackage{xcolor}
\usepackage[colorlinks,
  linkcolor=blue!60!black,
  citecolor=red!65!black,
  urlcolor=green!50!black]{hyperref}

\usepackage{bookmark}

\usepackage{mathtools}
\usepackage{amsthm}
\usepackage{amssymb}
\usepackage{mathrsfs}

\DeclarePairedDelimiter{\abs}{|}{|}
\DeclarePairedDelimiter{\floor}{\lfloor}{\rfloor}
\DeclarePairedDelimiter{\norm}{\Vert}{\Vert}
\DeclarePairedDelimiter{\Angle}{\langle}{\rangle}
\renewcommand{\angle}{\Angle}

\DeclarePairedDelimiter{\cat}{[}{]}
\newcommand{\E}{\mathbb{E}}
\renewcommand{\P}{\mathbb{P}}
\newcommand{\N}{\mathbb{N}}
\newcommand{\R}{\mathbb{R}}
\newcommand{\Z}{\mathbb{Z}}

\newcommand{\PP}{\mathbf{P}}
\newcommand{\QQ}{\mathbf{Q}}
\newcommand{\Qcal}{\mathcal{Q}}

\newcommand{\spine}{\mathbf{s}}
\newcommand{\indic}{\mathbf{1}}
\newcommand{\an}[1]{\overleftarrow{#1}}

\newcommand{\diff}{\mathrm{d}}

\DeclareMathOperator{\Leb}{Leb}
\DeclareMathOperator{\Card}{Card}

\newcommand{\goesto}[2][]{\xrightarrow[\:#2\:]{\:#1\:}}

\renewcommand{\phi}{\varphi}
\renewcommand{\epsilon}{\varepsilon}
\renewcommand{\emptyset}{\varnothing}

\newcommand{\tiff}{if\,\!f }

\allowdisplaybreaks

\newtheorem{theorem}{Theorem}
\newtheorem{proposition}{Proposition}
\newtheorem{corollary}{Corollary}
\newtheorem{lemma}{Lemma}

\theoremstyle{definition}
\newtheorem{definition}{Definition}
\newtheorem{assumption}{Assumption}
\newtheorem{remark}{Remark}

\title{Convergence of genealogies through spinal decomposition with an
application to population genetics}

\date{}

\usepackage{authblk}

\newlength{\affilskip}
\author[1]{Félix Foutel-Rodier}
\author[2]{Emmanuel Schertzer}

\affil[1]{
    Department of Statistics, University of Oxford, UK
    \vspace{\affilskip}
}
\affil[2]{
    Faculty of Mathematics, University of Vienna, \authorcr
    Oskar-Morgenstern-Platz 1, 1090 Wien, Austria
    \vspace{\affilskip}
}

\usepackage[hang,flushmargin]{footmisc} 
\makeatletter
\def\blfootnote{\xdef\@thefnmark{}\@footnotetext}
\makeatother

\begin{document}

\maketitle

\begin{abstract}
Consider a branching Markov process with values in some general type
space. Conditional on survival up to generation $N$, the genealogy of the
extant population defines a random marked metric measure space, where
individuals are marked by their type and pairwise distances are measured
by the time to the most recent common ancestor. In the present
manuscript, we devise a general method of moments to prove
convergence of such genealogies in the Gromov-weak topology when $N \to
\infty$. 

Informally, the moment of order $k$ of the population is obtained by
observing the genealogy of $k$ individuals chosen uniformly at random
after size-biasing the population at time $N$ by its $k$-th factorial
moment. We show that the sampled genealogy can be expressed in terms of a
$k$-spine decomposition of the original branching process, and that
convergence reduces to the convergence of the underlying $k$-spines.  

As an illustration of our framework, we analyse  the large-time behavior
of a branching approximation of the biparental Wright--Fisher model with
recombination. The model exhibits some interesting mathematical
features. It starts in a supercritical state but is naturally driven to
criticality.  We show that the limiting behavior exhibits both critical
and supercritical characteristics. 
\end{abstract}

\section{Introduction}

\blfootnote{\textbf{Keywords and phrases:} spinal decomposition;
many-to-few formula; convergence of genealogies; Yaglom's law; recombination;
population genetics\\
\textbf{AMS 2010 Classification:} Primary 60J80; Secondary 60G57, 60J85,
92D25}

\paragraph{A branching process with a self-organized criticality behavior.} 
The original motivation of the present article is the branching
approximation of a classical model in population genetics. It can 
be formulated as a branching process in discrete time where each
individual carries a subinterval of $(0, R)$, for some fixed parameter $R
> 0$. At generation $t = 0$, the population is made of a single
individual carrying the full interval $(0, R)$. At each subsequent
generation, individuals reproduce independently and an individual
carrying an interval $I$ with length $\abs{I}$ gives birth to $K(I)$
children, where
\[
    K(I) \sim \mathrm{Poisson}\big(1 + \tfrac{\abs{I}}{N}\big),
\]
and $N \ge R$ is another fixed parameter. Each of these $K(I)$ children
inherits independently an interval which is either the full parental
interval $I$, or a fragmented version of it. More precisely, with
probability
\[
    r_N(I) = 2 \frac{\abs{I}}{N} \big(1 + o_N(1)\big)
\]
we say that a \emph{recombination} occurs: a random point is sampled
uniformly on $I$ which breaks $I$ into two subintervals. The
child inherits either the left or the right subinterval with equal
probability. With probability $1 - r_N$ no recombination occurs and the
child inherits the full parental interval $I$. We refer to this process
as the \emph{branching process with recombination}.

One of the most interesting aspect of the present model is a
\emph{self-organized criticality} property. While the process is
``locally'' supercritical, since $\E[K(I)] > 1$, intervals are broken
via recombination and the process is naturally driven to criticality.
Under the regime $N \gg R \gg 1$, we will prove that some features are
reminiscent of a critical branching process (for instance, it satisfies a
type of Yaglom's law) but also bears similarities to supercritical
branching processes. In particular, one striking feature is related to
the genealogy of the process conditioned on survival at a large
time horizon. In the natural time scale, the genealogy of the extant
population is indistinguishable from the supercritical case, that is, it
converges to a star tree. However, if we zoom in on the root by rescaling
time in a logarithmic way, the genealogy converges to the celebrated
Brownian Coalescent Point Process and becomes indistinguishable from a
critical branching process.

From a biological standpoint, our process was first introduced in 
\cite{baird_distribution_2003} and corresponds to a branching
approximation of a more complicated model of population genetics, named
the biparental Wright--Fisher model with recombination. The connection
between the two models and their biological significance are discussed in
greater details in Section~\ref{S:maintResults}.

\paragraph{Convergence of types and genealogy.}
In order to analyse the previous model, we introduce a general framework
and provide simple criteria for the convergence of random genealogies. 
Although the branching process that we consider is interesting in its own
right, our study aims at giving a concrete illustration of a general
approach that could presumably be relevant in many other settings.

It is quite common that individuals in a branching process are endowed
with a ``type'', which is heritable and can in turn influence the
reproductive success of individuals. Let us denote by $E$ the set of
types. For instance, in our work $E$ is the set of subintervals of
$(0,R)$, for branching random walks $E = \R^d$ \cite{shi_2015}, or for
multi-type Galton--Watson processes $E$ is often chosen as finite
\cite[Chapter~5]{athreya1972}. In the absence of types or
when the reproduction law does not depend on types (as for standard
branching random walks in $\R^d$), the scaling limits of the tree
structure and of the distribution of types have received quite a lot of
attention \cite{duquesne2002, popovic2004, legall91}. In this particular
setting, one can make use of an encoding of the tree as the excursion of
a stochastic process, the so-called contour process, or height process.
Convergence is then obtained by showing that the corresponding excursion
converges.

When the reproduction law may depend on the types, some attempts to
extend the excursion approach exist in the literature
\cite{popovic2014coalescent} but as far we know a systematic and amenable
approach is still missing. In this work we follow a different approach,
and extend the seminal work of \cite{greven_convergence_2006} to prove
convergence in the Gromov-weak topology.
Proving convergence in distribution for
this setting is very similar in spirit to the method of moments for real
random variables, where one proves convergence in distribution by showing
that \emph{all moments of the tree structure} converge. In the context of
trees and metric spaces, the moments of order $k$ are obtained by summing
over all $k$-tuples of individuals at some generation, and considering a
functional of the subtree spanned by these $k$ individuals. Informally,
this amounts to picking $k$ individuals at random in a size biased
population, and then proving convergence of the genealogy of the sample.
One contribution of our work is that, analogously to the method of
moments in the real setting, we only need to prove convergence of the
moments with no need to identify the limit. This relies on a 
de~Finetti-like representation of exchangeable coalescents that was
developed in \cite{foutelrodier2019exchangeable}. See
Theorem~\ref{thm:convUMS} for our main convergence result.

\paragraph{Spinal decomposition of Markov branching processes.}
To compute the moments of branching process, we make use of a second set
of tools called spinal decompositions \cite{lyons_conceptual_1995,
shi_2015, harris_many_2017}. One of the main insight of the present
manuscript lies in the observation that an ingenious random change of
measure allows us to reduce the computation of a polynomial of order $k$
to a computation on a single tree with $k$ leaves, called the $k$-spine
tree. Since this type of manipulation allows one to reduce a computation
involving the whole tree to a computation involving only $k$ individuals,
this type of results have been called \emph{many-to-few formula}. While
\emph{many-to-one} formula have been extensively explored in the
literature since the seminal work of Lyons, Pemantle, and Peres
\cite{lyons_conceptual_1995}, spinal decompositions of higher order are
more sparse. One formulation has been exposed in \cite{harris_many_2017}
(see also \cite{johnston2019, harris2020, bansaye11, ren2018})
where the $k$-spine is constructed from a system of branching particles
evolving according to a prescribed Markov dynamics. While the main result
in \cite{harris_many_2017} could be in principle applied to our setting,
the computation rapidly proved to be intractable. Another contribution of
our work, that we want to emphasize, is a derivation of new general
many-to-few formula that was better suited to our case. Let us describe
it briefly here, and refer to Section~\ref{S:spinalDecomposition} for a
complete account. 

The $k$-spine tree in our work is constructed iteratively as a
\emph{coalescent point process} (CPP in short).  Starting from a single
branch of length $N$, at each step a new branch is added to the right of
the tree. Branch lengths are assumed to be i.i.d., and the procedure is
stopped when the tree has $k$ leaves, see Figure~\ref{fig:CPP}. Given
this tree, types need to be assigned to vertices of the $k$-spine tree.
For $k = 1$, the tree is made of a single branch, and the sequence of
types observed from the root to the unique leaf is a Markov chain. This
Markov chain is the usual sequence of types along the spine that arises
in many versions of the many-to-one formula \cite{biggins2004measure,
shi_2015}. It is obtained as the Doob harmonic transform of the offspring
type, see Section~\ref{SS:manyToFew}. For a general $k$, the previous
chain is duplicated independently at each branch point. The distribution
of the resulting tree is connected to the original distribution of the
branching process through a random change of measure $\Delta_k$ given in
\eqref{eq:bias}. The latter factor accounts for the fact that individuals
located at the branch points are more likely to have a large offspring
and a favorable type.

While our spinal decomposition result bears similarities with that in
\cite{harris_many_2017}, our formulation allows for a more general
distribution of the $k$-spine tree, which can be any discrete CPP.
This additional degree of freedom proved very valuable in our
application, where the introduction of a well-chosen ansatz for the
genealogy of the process, see \eqref{eq:branchTimes}, simplified
considerably earlier versions of our proofs. More generally, we believe
that our approach is particularly amenable to the study of near-critical
branching processes, since the scaling limit of their genealogy can also
be described as a continuous CPP. Nevertheless, see \cite{johnston2019,
harris2020} for successful applications of the techniques in
\cite{harris_many_2017} to study the genealogy of a sample from a
Galton--Watson tree.


\paragraph{Outline.}
Overall, the contribution of our work is three-fold. We have 1) derived a
new type of many-to-few formula based on a CPP tree, 2) combined it with
the framework of the Gromov-weak topology to produce an effective way of
studying the scaling limit of types and genealogies in branching
processes, and 3) applied it to study a complex model from population
genetics, the branching process with recombination.

The rest of our work is laid out as follows. Section~\ref{S:maintResults}
provides more details on the biological motivation of the branching
process with recombination and a statement of our main results concerning
this model. Those results will be proved using a general framework that
will be developed in the subsequent sections.

In Section~\ref{S:spinalDecomposition} we construct the $k$-spine tree
and prove our spinal decomposition result. In Section~\ref{S:GromovWeak},
we show that the convergence of the genealogy of branching processes
can be reduced to the convergence of the associated $k$-spines. This
approach relies on a previous work \cite{foutelrodier2019exchangeable}
where we provide a de~Finetti-like representation of ultrametric spaces
that allows us to extend previous convergence criteria for the
Gromov-weak topology.

In the last two sections, we apply the previous framework to the model at
hand.  In Section~\ref{S:1spine}, we characterize the $1$-spine
associated to the branching process with recombination, and prove our
convergence results in Section~\ref{S:kspine}.

\section{Branching process with recombination}
\label{S:maintResults}

\subsection{Biological motivation}

\begin{figure}
    \centering
    \includegraphics[width=.9\textwidth]{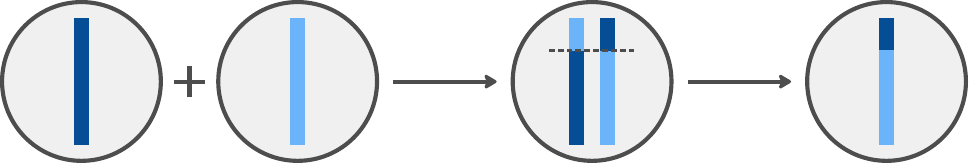}
    \caption{When a recombination occurs, a point is chosen along the
    sequence, called the crossover point and represented by a dashed
    line. Two new chromosomes are formed by swapping the parts of the
    parental chromosomes on one side of the crossover point. The offspring
    inherits one of these two chromosomes.}
    \label{fig:crossover}
\end{figure}

In the context of this work, genetic recombination is the biological
mechanism by which an individual can inherit a chromosome which is not a
copy of one of its two parental chromosomes, but a mix of them. An
idealized version of this mechanism is illustrated in
Figure~\ref{fig:crossover}. Due to recombination, the alleles carried by
an individual at different loci, that is, locations on the chromosome,
are not necessarily transmitted together. At the level of the population,
this creates a complex correlation between the gene frequencies at
different loci which is hard to study mathematically.

When focusing on a finite number of loci it is possible to express the
dynamics of these frequencies as a set of non-linear differential
equations or stochastic differential equations \cite{Baake2008,
ohta_kimura_1969}. However, one needs to keep track of the frequencies of
all possible combinations of alleles. As the number of such combinations
grows exponentially fast with the number of loci, it leads to expressions
that rapidly become cumbersome, providing little biological insight.
Another very fruitful approach is to trace backward-in-time the set of
potential ancestors of the population. This gives rise to a mathematical
object named the ancestral recombination graph (ARG)
\cite{griffiths1997}, or see \cite[Chapter~3]{durrett2008}. However, the
ARG is quite complicated both from a mathematical and a numerical point
of view. Nevertheless, see \cite{lambert21} for some recent mathematical
results, and \cite{mcvean2005, Marjoram2006} for approximations of the
ARG that have proved very successful in application.

In this work we consider a third approach to this question, which is to
envision the chromosome as a continuous segment. At each reproduction
event recombination can break this segment into several subintervals, a
subset of which is transmitted to the offspring, as in
Figure~\ref{fig:crossover}. The genetic contribution of an individual is
now described by a collection of intervals, which are delimited by points
called \emph{junctions}. This point of view has long-standing history
dating back to the work of Fisher \cite{Fisher1954}, see for instance
\cite{janzen21} and references therein. Let us discuss the specific model
that we consider, and how the branching process with recombination
approximates it.

\subsection{Connection to the Wright--Fisher model} 

\begin{figure}
    \centering
    \includegraphics[width=.9\textwidth]{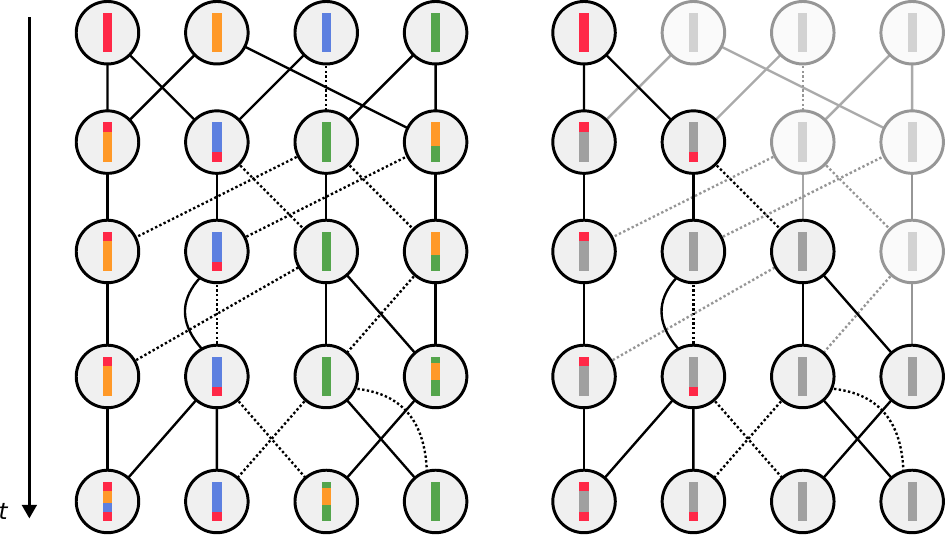}
    \caption{Illustration of the Wright--Fisher model with recombination.
    A line is drawn between each individual and its parents. It is dotted
    if no genetic material is inherited from this parent. The right panel
    focuses on the genetic material left by the red individual. Note that
    each individual only carries an interval before the first reproduction
    event involving two descendants of the focal ancestor.}
    \label{fig:wrightFisher}
\end{figure}

Consider a population of fixed size $N$ where individuals are endowed
with a continuous chromosome represented by the interval $(0, R)$. At
each generation, individuals pick independently two parents uniformly
from the previous generation. Assume that these parents can be
distinguished, so that there is a left and a right parent. Then,
independently for each individual:
\begin{itemize}
    \item with probability $1 - R/N$, it inherits the chromosome of one
        of its two parents, say the left one;
    \item with probability $R/N$, a recombination occurs. A crossover
        point $U$ is sampled uniformly on $(0, R)$, and the offspring
        inherits the part of the chromosome to the left of $U$ from its
        left parent, and that to the right of $U$ from its right parent.
\end{itemize}
Suppose that at some focal generation, labeled generation $t = 0$, each
chromosome in the population is assigned a different color. Due to
recombination, new chromosomes are formed that are mosaics of the initial
colors. We are ultimately interested in describing the long-term
distribution of these mosaics in the population. This is
illustrated in Figure~\ref{fig:wrightFisher}.

In this work, we consider a simpler but related problem. Fix a focal
ancestor, and say that its chromosome is red. We trace the individuals in
the population that have inherited some genetic material from this focal
ancestor, that is, the set of individuals that have some red on
their chromosome as well as the location of the red color.
To recover the branching approximation that we study, consider an
individual in the population at some generation $t$ carrying a red
interval $I$. Its offspring size distribution is
\[
    \mathrm{Binomial}\Big(N, \frac{2}{N}\big(1 + o_N(1)\big)\Big)
    \goesto{N \to \infty} \mathrm{Poisson}(2).
\]
Each of these children has another parent in the population. As long as
the number of individuals with a red piece of chromosome is small
compared to $N$, this other parent does not have any red part on its
chromosome.

Therefore, there are only four possible outcomes for each child:
\begin{itemize}
    \item With probability $1 - R/N$ no recombination occurs and
        \begin{itemize}
            \item with probability $1/2$ it inherits $I$; 
            \item with probability $1/2$ the interval $I$ is lost.
        \end{itemize}
    \item With probability $R/N$ a recombination occurs and
        \begin{itemize}
            \item if $U \not\in I$ the interval $I$ is transmitted or
                lost with probability $1/2$; 
            \item if $U \in I$, the child inherits the subinterval of $I$
                to the left or to the right of $U$ with probability
                $1/2$.
        \end{itemize}
\end{itemize}
By combining the previous cases, we recover that the number of
descendants that carry some red of an individual with red interval $I$
is approximately distributed as a $\mathrm{Poisson}\big(1 + \tfrac{\abs{I}}{N}
\big)$ r.v., and that the probability of inheriting a fragmented
interval is
\[
    r_N(I) = \frac{2 \abs{I}/N}{1+\abs{I}{N}}.
\]
This is the description of the branching process with recombination.

\subsection{Limiting behavior}
\label{SS:mainResults}

Let $\PP_R$ denote the distribution of the branching process with
recombination started from a single individual with interval $(0, R)$.
The following asymptotic expression for the survival probability of this
process was already derived in \cite{baird_distribution_2003}.

\begin{proposition}[\cite{baird_distribution_2003}] \label{prop:survival}
    Let $Z_N$ denote the population size at generation $N$ in the
    branching process with recombination. The limit 
    \[
        \lim_{N \to \infty} \frac{\PP_R(Z_N > 0)}{N}
    \]
    exists. Moreover it fulfills
    \[
        \lim_{R \to \infty} \lim_{N \to \infty} 
        \frac{R \, \PP_R(Z_N > 0)}{N\log R} = 1.
    \]
\end{proposition}

Let $T_N$ denote the set of individuals at generation $N$ in the
branching process with recombination. For an individual $u \in T_N$, we
denote by $I_u$ the interval that it carries. For $u, v \in T_N$, let
$d_T(u,v)$ denote the genealogical distance between $u$ and $v$, that is,
the number of generations that need to be traced backward in time before
$u$ and $v$ find a common ancestor.

Our first result provides the joint limit of the interval lengths and of
the genealogy of the population. To derive this limit, we will envision
the population as a marked metric measure space and work with the
marked Gromov-weak topology \cite{depperschmidt_marked_2011}. The
definition of this topology is recalled in
Section~\ref{SS:gromovWeak}.

Let us consider the measure $\mu_N$ on $T_N \times \R_+$ defined as 
\[
    \mu_N = \sum_{u \in T_N} \delta_{(u, \abs{I_u})}.
\]
The triple $[T_N, d_T, \mu_N]$ is the marked metric measure space
corresponding to the branching process with recombination. Let us finally
define the rescaling
\begin{equation} \label{eq:rescaling}
    \forall x \in [0, 1], \quad 
    F_R(x) = \frac{\log \big( (R-1)x + 1\big)}{\log R},
    \quad F^{-1}_R(x) = \frac{R^x-1}{R-1}
\end{equation}
and define the rescaled distance as 
\[
    \forall u,v \in T_N,\quad \bar{d}_N^R = 1 - F_R\big( 1 - \tfrac{d_T(u,v)}{N} \big)
\]
which is the distance obtained by rescaling time according to $F_R$.

\begin{theorem} \label{thm:main1}
    Fix $t > 0$. Conditional on survival at time $\floor{Nt}$ the
    following limit holds in distribution for the marked Gromov-weak
    topology,
    \[
        \lim_{R \to \infty} \lim_{N \to \infty} 
        \Big[ T_{\floor{Nt}}, \bar{d}_{\floor{Nt}}^R,
            \frac{\mu_{\floor{Nt}}}{t N \log R} \Big]
        = \big[ (0, Y), d_P, \Leb \otimes \mathrm{Exp}(t) \big]
    \]
    where $[(0,Y), d_P]$ is a Brownian coalescent point process, and
    $\mathrm{Exp}(t)$ is the exponential distribution with mean $1/t$. 
\end{theorem}

A stronger version of this result is proved in
Section~\ref{SS:mainProof}. Let us now briefly discuss several
consequences of the previous result.

\paragraph{Convergence of the empirical measure.} As
mentioned in the introduction, the branching process at hand is
naturally driven to criticality through recombination. Recall that if
the offspring distribution of a (standard) critical branching
processes has finite second moment, the celebrated Yaglom law
states that conditional on survival up to time $Nt$, the rescaled
population size $Z_{\floor{Nt}}/N$ converges to an exponential random
variable. In contrast, the convergence of $\frac{\mu_{\floor{Nt}}}{t
N \log R}$ entails that the rescaled population size at time $Nt$
converges to an exponential random variable, but the population size
is of order $N\log R$ instead of $N$. In words, the local
supercritical character of the process translates into an extra
$\log R$ factor for the population size.  
 
Secondly, the convergence of the random measure
$\frac{\mu_{\floor{Nt}}}{t N \log R}$ also implies that the length of the
interval carried by a typical individual in the population is
exponentially distributed with mean $1/t$. Since the limiting random
measure is deterministic, the intervals carried by $k$ typical
individuals in the populations are independent (\emph{propagation of
chaos}). Note that, although the length of the initial interval $R$ goes
to infinity, the intervals at any finite time $t$ remain of finite
length. This phenomenon is usually referred to as \emph{coming down from
infinity}. In our work, it originates from the existence of an entrance
law at infinity for the spine, which turns out to be connected to the
existence of such entrance laws for positive self-similar Markov
processes with negative index of self-similarity
\cite{bertoin02,bertoin2002}.

\paragraph{Convergence of the genealogy.} Let us first comment on the
rescaling $F_R$. Although the expression of $F_R$ appears a bit daunting
at first, it essentially boils down to first rescaling time  by $N$ (as
expected), and then  measuring time from the origin in the log-scale. The
first consequence is that the genealogy of the population in the natural
scale (that is, if we only rescale time by $N$) converges to a star tree so
that the genealogy becomes indistinguishable from the one of a
supercritical branching process at the limit.

\begin{figure}[t]
    \centering
    \includegraphics[width=.8\textwidth]{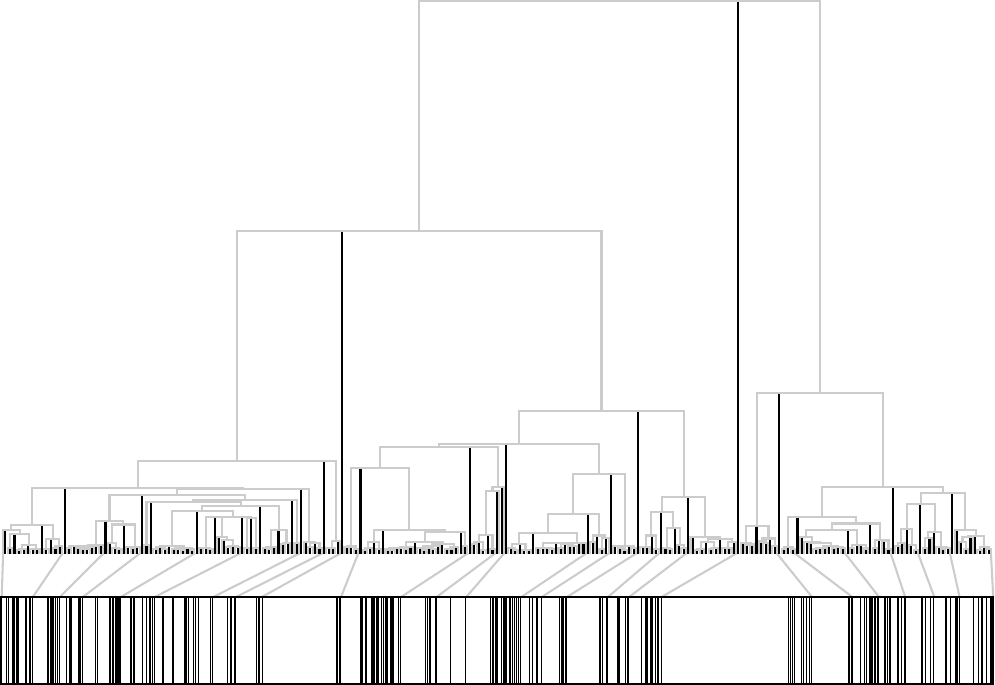}
    \caption{Top: simulation of a Brownian CPP. The black vertical lines represent
        to the atoms of $P$, and the corresponding tree is pictured
        in grey. Bottom: geometry of the blocks of ancestral material corresponding
        to the top CPP. Each block is represented by a black stripe. The
        correspondence between the blocks and the tree are shown for some
        blocks by grey segments joining the two. The distance between two
        consecutive stripes is the logarithm of their distance on the
        chromosome. Note that this induces a strong deformation of the
        intuitive linear scale.}
    \label{fig:geometry}
\end{figure}

A second consequence of this result is that, after rescaling time
according to $F_R$, the genealogy of the branching process with
recombination converges to a limiting metric space named the Brownian
coalescent point process (CPP). It is constructed out of a Poisson point
process $P$ on $(0,\infty) \times (0, 1)$ with intensity $\diff t \otimes
\frac{1}{x^2} \diff x$. Let 
\[
    \forall x \le y,\quad d_P(x,y) = \sup \{ z : (t,z) \in P,\, x \le t \le y \}.
\]
The Brownian CPP is the random metric space $[(0,Y), d_P]$, where $Y$ is
an exponential r.v.\ with mean $1$, independent of $P$, see
Figure~\ref{fig:geometry} for a graphical construction. It corresponds to
the limit of the genealogy of a critical Galton--Watson process with
finite variance \cite{popovic2004}.

\paragraph{Chromosomic distance.} 
The previous result provides a complete description of the interval
lengths in the population, but does not provide any insight into their
distribution over $(0, R)$. We will encode the latter information by
picking a reference point belonging to each interval in the population
and considering the usual distance on the real line between these points.
More precisely, for each $u \in T_N$, pick a reference point $M_u$
uniformly on $I_u$. We define a new metric 
\[
    \forall u,v \in T_N,\qquad D_N(u,v) = \abs{M_u - M_v}.
\]
We will refer to the $D_N$ as the \emph{chromosomic distance}. 

The quadruple $[T_N, d_N, D_N, \mu_N]$ can be seen as a random
``bi-metric'' measure space with marks. We can define a straightforward
extension of the marked Gromov-weak topology for such objects, see the
end of Section~\ref{SS:gromovWeak}. The correct rescaling for $D_N$ is to
set 
\[
    \forall u,v \in T_N, \quad \bar{D}_N^R = \frac{\log D_N \vee 2}{\log R} 
\]
In Section~\ref{SS:mainProof}, we prove the following refinement of
Theorem~\ref{thm:main1}.

\begin{theorem} \label{thm:main2}
    Fix $t > 0$. Conditional on survival of the process at $\floor{tN}$, 
    the following limit holds in distribution for the marked Gromov-weak
    topology,
    \[
        \lim_{R \to \infty} \lim_{N \to \infty} 
        \Big[ T_{\floor{Nt}}, \bar{d}_{\floor{Nt}}^R, \bar{D}_{\floor{Nt}}^R, 
            \frac{\mu_{\floor{tN}}}{t N \log R} \Big]
        = \big[ (0, Y), d_P, d_P, \Leb \otimes \mathrm{Exp}(t) \big]
    \]
    where $[(0,Y), d_P]$ is a Brownian coalescent point process, and
    $\mathrm{Exp}(t)$ is the exponential distribution with mean $1/t$. 
\end{theorem}

It is important to note that, in the limit, the two metrics coincide.
This result is quite interesting from a biological point of view. It
shows that there is a correspondence between the genealogical distance
between two individuals, and the chromosomic distance between the genetic
material that they carry. Indeed, the latter two quantities are
correlated: two individuals inherit intervals that are subsets of the
interval carried by their most-recent common ancestor. If this ancestor
is recent, its interval is smaller, and so is their chromosomic distance.
Our result shows that, in the limit, the two distances become identical
when considered on the right scale. This result is illustrated in
Figure~\ref{fig:geometry}.

\begin{remark}
    Define the point process
    \[
        \vartheta = \sum_{x \ge 0} \Leb(\{ y \in (0, Y) : d_P(0, y) = x \}) \delta_x
    \]
    which corresponds to the CPP tree ``viewed from the individual with
    the left-most interval''. Using elementary properties of Poisson
    point processes shows that $\vartheta$ can also be written as 
    \[
        \vartheta = \sum_{(x_i, y_i) \in \mathcal{P}} y_i \delta_{x_i}
    \]
    where $\mathcal{P}$ is a Poisson point process on $(0,\infty) \times
    (0, \infty)$ with intensity $\frac{1}{x^2} e^{-x/y} \diff x \diff y$.

    The same expression was obtained in \cite[Theorem~1.5]{lambert21} to
    describe the set of loci that share the same ancestor as the
    left-most locus in the fixed haplotype of a Wright--Fisher model with
    recombination, under a limiting regime similar to ours. This
    connection is quite surprising. We are considering a branching
    approximation where all intervals belong to distinct individuals
    and its chromosome carries at most one block of ancestral genome,
    whereas in \cite{lambert21} all intervals belong to a single
    chromosome, which has reached fixation in the population.
\end{remark}

\section{The \texorpdfstring{$k$}{k}-spine tree}
\label{S:spinalDecomposition}

\subsection{The many-to-few formula}
\label{SS:manyToFew}

The objective of this first section is to introduce the $k$-spine tree
and state our many-to-few formula, that relates the expression of the
polynomials of a branching process to the $k$-spine tree. All the random
variables introduced here are more formally defined in the forthcoming
sections, where the proof of the many-to-few formula is carried out. A
formal statement of our result requires some preliminary notation.
 
\paragraph{Assumption and notation.} Consider a Polish space $(E, d_E)$,
and a collection $(\Xi(x);\, x \in E)$ of random point measures on $E$.
This collection can be used to construct a branching process with type
space $E$, such that the atoms of a realization of $\Xi(x)$ provide the
types of the children of an individual with type $x$. The distribution of
the resulting branching process is denoted by $\PP_x$.

Let $K(x)$ denote the number of atoms of $\Xi(x)$, and set
\[
    \pi(x,k) = \P(K(x) = k)
\]
for the distribution of $K(x)$. The $n$-th factorial moment of $K(x)$ is
denoted by $m_{n}(x)$, that is,
\begin{gather*}
    m_{k}(x) \coloneqq \E\Big[ K(x)^{(n)} \Big],\\
    m(x) \coloneqq m_{1}(x) = \E[K(x)],
\end{gather*}
where we have used the notation $k^{(n)}$ for the $n$-th descending
factorial of $k$, 
\[
    k^{(n)} = k(k-1) \dots (k-n+1).
\]
Our results are more easily formulated under the assumption that,
conditional on $K(x)$, the locations of the atoms are i.i.d.\ with
distribution $p(x, \cdot)$. That is, we assume that
\[
    \Xi(x) = \sum_{i = 1}^{K(x)} \delta_{\xi_i(x)}
\]
where $(\xi_i(x);\, i \ge 1)$ is an i.i.d.\ sequence distributed as
$p(x, \cdot)$ and is independent of $K(x)$. We make the further
simplifying assumption that all distributions $p(x, \cdot)$ have a
density w.r.t.\ some common measure $\Lambda$ on $E$. With a slight abuse
of notation, the density of $p(x,\cdot)$ is denoted by $(p(x,y);\, y \in E)$.

\paragraph{Harmonic function.} We say that a map $h \colon E \to
[0,\infty)$ is (positive) harmonic if 
\[
    \forall x \in E,\quad h(x) = \E[\Angle{\Xi(x), h}],
\]
where we used the notation $\Angle{\mu, f} = \int f \diff \mu$, see for
instance \cite{biggins2004measure}. A harmonic function can be used to
define a new probability kernel on $E$, defined as 
\begin{equation} \label{eq:harmonicKernel}
    \forall x,y \in E,\quad q(x,y) = \frac{m(x) h(y) p(x,y)}{h(x)}.
\end{equation}
The fact that this is a probability measure follows from the
harmonicity of $h$.

\begin{figure}
    \centering 
    \includegraphics[width=.5\textwidth]{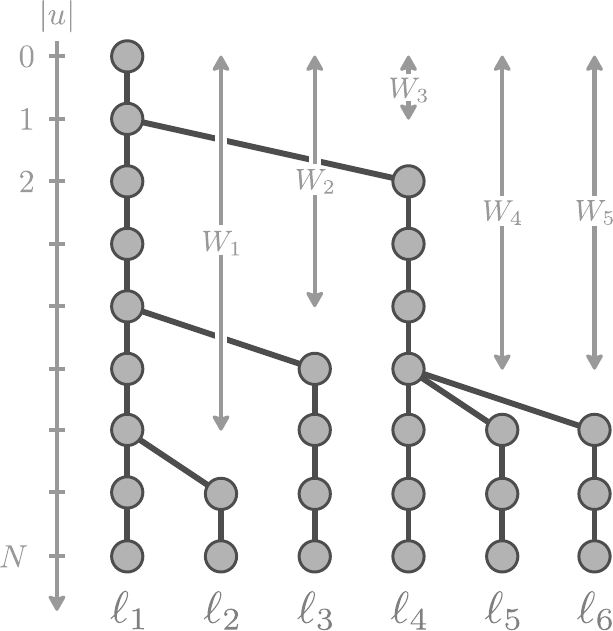}
    \caption{Illustration of the construction of a CPP tree. The vector
        $(W_1,\dots,W_{k-1})$ branching times between successive leaves.
        In this example, this vector is $(6,4,1,5,5)$. The tree is
        recovered from these times by grafting for each $i$ a branch of
        length $N-W_i$ to the right-most vertex of the tree at generation
        $W_i$.}
    \label{fig:CPP}
\end{figure}

\paragraph{The $k$-spine tree.} We are now ready to define the $k$-spine
tree. Let $\nu = (\nu_n;\, n \in \{0, \dots, N-1\})$ be a probability
distribution and let $(W_1, \dots, W_{k-1})$ be i.i.d.\ random variables 
with distribution $\nu$. Define 
\[
    \forall i \le j,\quad d_T(i,j) = d_T(j,i) = \max \{ N-W_i, \dots, N-W_{j-1} \}.
\]
There is a unique tree with $k$ leaves labeled by $\{1,\dots,k\}$ such
that the tree distance between the leaves is $d_T$. We denote it by $S$
and call it the $\nu$-CPP tree. This tree is constructed inductively by
grafting a branch of length $N-W_i$ on the tree constructed at step $i$,
as illustrated in Figure~\ref{fig:CPP}.

We now assign marks on the tree such that along each branch of the
tree, marks evolve according to a Markov chain with transition kernel
$(q(x,y);\, x,y \in E)$ defined in \eqref{eq:harmonicKernel}. More
formally construct a collection of processes $(X_1, \dots, X_k)$ such
that 
\begin{itemize}
    \item the process $(X_1(n);\, n \ge 0)$ is a Markov chain with
        transition $(q(x,y);\, x,y \in E)$ started from $x$;
    \item conditional on $(X_1, \dots, X_i)$, 
        \[
            \forall n \ge 0,\quad 
            X_{i+1}(n) =
            \begin{cases}
                X_i(n) &\text{if $n < W_i$}\\
                X'(n-W_i) &\text{if $n \ge W_i$}
            \end{cases}
        \]
        for some independent Markov chain $X'$ with transition 
        $(q(x,y);\, x,y \in E)$ started from $X_i(W_i)$.
\end{itemize}
By thinking of $(X_i(n);\, n \ge 0)$ as giving the sequence of marks
along the branch of $S$ starting from the root and going to the $i$-th
leaf, we can assign to each vertex $u \in S$ a mark $Y_u$.

\begin{definition}
    The $k$-spine tree is the random marked tree $[S, (Y_u;\, u \in S)]$
    encoded by the variables $(W_1, \dots, W_{k-1})$ and $(X_1, \dots,
    X_k)$. The distribution of the latter variables is denoted by
    $\QQ^{k,N}_x$.
\end{definition}

We are now ready to state our many-to-few formula. It can be described
informally as follows. Suppose that the branching process with law
$\PP_x$ is biased by the $k$-th factorial moment of its size at
generation $N$ and that $k$ individuals are chosen uniformly for that
generation. Then the law of the subtree spanned by these individuals is
$\QQ^{k,N}_x$ biased by a random factor $\Delta_k$ that can be expressed
as 
\begin{equation} \label{eq:bias}
    \Delta_k = \prod_{\substack{u \in S\\ d_u > 1}} 
        \Big(\frac{ h(Y_u)}{N\nu_{\abs{u}}} \Big)^{d_u-1} \frac{m_{d_u}(Y_u)}{d_u!\,m(Y_u)^{d_u}}
            \cdot \prod_{i=1}^k \frac{1}{h(X_i(N))},
\end{equation}
where $d_u$ denotes the degree of a vertex $u \in S$ and $Y_u$ its mark.
Note that the left product in \eqref{eq:bias} has at most $k-1$ terms,
which correspond to the branch points in $S$.

\begin{lemma}\label{lem:simplify}
    Assume that for every $x\in E$, the offspring number $K(x)$ is Poisson
    (for some given parameter $\lambda(x)>0$ that may depend on $x$).  Then 
    \[
        \Delta_k = \prod_{\substack{u \in S\\ d_u > 1}} 
            \Big(\frac{h(Y_u)}{N\nu_{\abs{u}}} \Big)^{d_u-1} \frac{1}{d_u!}
            \cdot \prod_{i=1}^k \frac{1}{h(X_i(N))}.
    \]
\end{lemma}

\begin{proof}
    This simply follows from the well known fact that the $k$-th
    factorial moment of a Poisson random variable with parameter
    $\lambda>0$ is $\lambda^k$.
\end{proof}

Finally, let $T_N$ denote the labels of the $N$-th generation of a
branching process with distribution $\PP_x$, for $u \in T_N$ let $X_u$
denote its type, and let $d_T$ denote the tree distance on $T_N$.

\begin{proposition}[Many-to-few] \label{prop:manyToFew}
    For any test function $\phi$, 
    \begin{multline*}
        \PP_x\Big[ 
            \sum_{\substack{(v_1,\dots,v_k) \in T_N\\ v_1 \ne \dots \ne v_k}} 
            \phi\big( d_T(v_i,v_j), X_{v_i};\, i,j \le k \big)
        \Big] \\
        = h(x) N^{k-1} k!\, \QQ^{k,N}_x\Big[ 
            \Delta_k \cdot 
            \phi\big( d_T(\sigma_i, \sigma_j), X_{\sigma_i}(N);\, i,j \le k\big)
        \Big]
    \end{multline*}
    where $\sigma$ is an independent uniform permutation of $\{1, \dots, k\}$.
\end{proposition}

\begin{remark}
    \begin{enumerate}[(i)]
        \item 
    In our construction of the $k$-spine, the distribution of the tree is
    independent of the marking. The term $\Delta_k$ captures the interplay
    between the genealogy and the types as a function of the marking at
    ``topological'' points. 

        \item 
    Compare Proposition~\ref{prop:manyToFew} to the many-to-few formula
    in \cite{harris_many_2017}. Both expressions relate the distribution
    of a $k$-sample from the branching process (l.h.s.\ of the equality) to that
    of a simpler $k$-spine tree (r.h.s.\ of the equality) at the expense of a
    bias term, here denoted by $\Delta_k$. 

        \item 
    In \cite{harris_many_2017}, the $k$-spine tree only depends on the
    moments of the reproduction law. Our formulation has one extra degree
    of freedom, since the $k$-spine tree is constructed out of an a
    priori genealogy, the $\nu$-CPP tree.

        \item 
    In many situations, including the model at hand, the bias term in
    \cite{harris_many_2017} becomes degenerate in the limit so that the
    distribution of the limiting genealogy is singular with respect to
    that of the original $k$-spine tree. For instance, for near-critical
    processes conditioned on survival at generation $N$, the first split
    time of the $k$-spine tree in \cite[Section~8]{harris_many_2017}
    remains of order $1$, whereas the most-recent common ancestor of the
    whole population is known to live at a time of order $N$. In
    contrast, one advantage of our approach is that $\nu$ can be
    well-chosen so that the bias $\Delta_k$ converges to a non-degenerate
    limit. This amounts to finding a good ansatz for the limiting
    genealogy. In our example this ansatz is given in
    \eqref{eq:branchTimes}, and the limit of the bias $\Delta_k$ is
    independent of the genealogy. This indicates that the limit of the
    genealogy does not depend on the types in the population.
    \end{enumerate}
\end{remark}

The rest of the section is dedicated to the proof of the many-to-few
formula. Our strategy to prove this result is to define a new tree with
distribution $\bar{\QQ}^{k,N}_x$ by grafting on the $k$-spine tree
independent subtrees distributed as the original branching process. The
many-to-few formula will then follow from the more precise spinal
decomposition theorem, which states that $\bar{\QQ}^{k,N}_x$ and $\PP_x$
are connected through the random change of measure $\Delta_k$. It is proved in
Section~\ref{SS:spinalDecomposition}. The remaining sections provide a
rigorous construction of the measures $\PP_x$ and $\bar{\QQ}_x^{k,N}$ and
the proof of the spinal decomposition theorem.

\subsection{Tree construction of the branching process}
\label{SS:notation}

Let us recall some common notation on trees.

\paragraph{Trees.} Following the usual Ulam--Harris labeling convention,
all trees will be encoded as subsets of 
\[
    \mathscr{U} \coloneqq \{ \emptyset \} \cup \bigcup_{n \ge 1} \N^n.
\]
Let us consider an element $u = (u(1), \dots, u(n)) \in \mathscr{U}$. We
denote by $\abs{u} = n$ its length, interpreted as the generation of $u$.
Moreover, its $i$-th child is denoted by
\[
    ui \coloneqq (u,i) = (u(1),\dots, u(n), i),
\]
and its ancestor in the previous generation as
\[
    \an{u} \coloneqq (u(1), \dots, u(n-1)).
\]
The set $\mathscr{U}$ is naturally endowed with a partial order
$\preceq$, where $u \preceq v$ if $u$ is an ancestor of $v$, that is,
\[
    u \preceq v \iff \forall i \le \abs{u},\; u(i) = v(i).
\]
The most-recent common ancestor of $u$ and $v$ can then be defined as 
\[
    u \wedge v \coloneqq \max \{w : w \preceq u \text{ and } w \preceq v \}.
\]
In the tree interpretation of $\mathscr{U}$, we can define a metric 
$d_T$ corresponding to the graph distance as 
\[
    \forall u,v \in \mathscr{U},\quad d_T(u,v) \coloneqq 
    \Card \{ w : u \wedge v \preceq w \prec u \}
    + 
    \Card \{ w : u \wedge v \preceq w \prec v \}.
\]
Finally, as a consequence of the Ulam--Harris encoding, trees are
\emph{planar} in the sense that the children of each vertex are endowed with
a total order. Accordingly let us denote by $\le$ the lexicographical
order on $\mathscr{U}$, which we will call the planar order. Note that 
$\le$ extends $\preceq$.

A subset $\tau \subseteq \mathscr{U}$ is called a \emph{tree} if
\begin{enumerate}[(i)] 
    \item $\emptyset \in \tau$;
    \item if, for some $j$, $uj \in \tau$, then $u \in \tau$;
    \item for any $u \in \tau$, there exists $k_u \in \N$ such that 
        \[
            ui \in \tau \iff i \le k_u,
        \]
        where $k_u$ is the number of children of $u$, also
        called the (out-)degree of $u$.
\end{enumerate}
The set of all trees is denoted by $\Omega$. For a tree $\tau \in
\Omega$, define its restriction to the $n$-th generation as
\[
    \tau_n = \{ u \in \tau : \abs{u} = n \}
\]
and that to the first $n$ generations as 
\[
    \tau_{[n]} = \{ u \in \tau : \abs{u} \le n \} = \bigcup_{i=1}^n \tau_i.
\]
Furthermore, let us denote by $\Omega_n$ the set of trees of height at
most $n$, where the height of a tree is defined as the generation of the
oldest individual in the tree.

\paragraph{Marked trees and definition of $\PP_x$.} A \emph{marked tree}
is a tree $\tau \in \Omega$ with a collection $(x_u;\, u \in \tau)$ of
marks with values in $E$. Let us define a random marked tree 
$[T, (X_u;\, u \in T)]$ inductively as follows, that corresponds to the
branching process with offspring reproduction point processes 
$(\Xi(x);\, x \in E)$.

Start from a single individual $\emptyset$ with mark $X_\emptyset = x$.
Conditional on the first $n$ generations $T_{[n]}$ and their marks $(X_u;\,
u \in T_{[n]})$, consider a collection of independent point processes
$(\Xi_u;\, u \in T_n)$, where 
\[
    \Xi_u \sim \Xi(X_u).
\]
Let us write 
\[
    \Xi_u = \sum_{i = 1}^{K_u} \delta_{\xi_{ui}}
\]
for the atoms of $\Xi_u$. Then define the next generation as
\[
    T_{[n+1]} = T_{[n]} \cup \{ ui : u \in T_n, i \le K_u \},
\]
with marks given by
\[
    \forall ui \in T_{n+1},\quad X_{ui} = \xi_{ui}.
\]
Let $T = \cup_{n \ge 1} T_n$ be the whole tree, and define define $\PP_x$
as the law of the random marked tree $[T, (X_u;\, u \in T)]$ obtained
through the previous procedure, and $\PP^N_x$ the law of its restriction
to the first $N$ generations.

\subsection{Ultrametric trees}

From now on, we consider a fixed, focal generation $N$. In this section
we construct the measure $\bar{\QQ}^{k,N}_x$ obtained by grafting some
independent subtrees on the $k$-spine tree. This construction relies on
the notion of (discrete) ultrametric trees.

\paragraph{Ultrametric trees.} A tree $\tau \in \Omega_N$ with height $N$
is called \emph{ultrametric} if all of its leaves lie at height $N$,
that is 
\[
    \forall u \in \tau, \quad k_u = 0 \implies \abs{u} = N.
\]
The set of all ultrametric trees of height $N$ with $k$ leaves is denoted
by $\mathbb{U}^{k,N}$. For $\tau \in \mathbb{U}^{k,N}$, let us denote by
$(\ell_1, \dots, \ell_k)$ the leaves of $\tau$ in lexicographical order,
that is, such that 
\[
    \ell_1 \le \dots \le \ell_k.
\]

The previous description of an ultrametric tree as an element of $\Omega_N$
is not suitable to describe the large $N$ limit of the $k$-spine tree.
To derive such a limit, we need to encode elements of $\mathbb{U}^{k,N}$
as a sequence $(g_1,\dots, g_{k-1})$ giving the branch times
between successive leaves in the tree. This construction is sometimes
referred to as a coalescent point process (CPP) \cite{popovic2004,
lambert13}.

More precisely, define the map
\[
    \Phi \colon \begin{cases}
        \mathbb{U}^{k,N} \to \{0, \dots, N-1\}^{k-1} \\
        \tau \mapsto \big( \abs{\ell_1 \wedge \ell_2}, \dots, \abs{\ell_{k-1} \wedge \ell_k} \big).
    \end{cases}
\]
The following straightforward result shows that the tree $\tau$ can be
recovered from the vector of coalescence times $\Phi(\tau)$.

\begin{lemma}
    The map $\Phi$ is a bijection from the set of ultrametric trees 
    $\mathbb{U}^{k,N}$ to the set of vectors $\{0, \dots, N-1\}^{k-1}$.
\end{lemma}

\paragraph{The $k$-spine tree.}
Let $(W_1,\dots, W_{k-1})$ and $(X_1,\dots,X_k)$ have distribution
$\QQ^{k,N}_x$. We formally define the $\nu$-CPP tree illustrated in
Figure~\ref{fig:CPP} as the random tree $S \coloneqq \Phi^{-1}(W_1,
\dots, W_{k-1})$. Note that the CPP tree associated to the uniform
distribution is uniform on $\mathbb{U}^{k,N}$. The processes
$(X_1,\dots,X_k)$ can now be used to construct a collection of marks
$(Y_u;\, u \in S)$ as follows. Each $u \in S$ is of the form
\[
    u = (\ell_i(1), \dots, \ell_i(n))
\]
for some leaf $\ell_i$ and $n \le N$. Define the mark of such a $u$ as 
\[
    Y_u \coloneqq X_i(n).
\]
(It is not hard to see that $Y_u$ is well-defined in that it does not
depend on the choice of $\ell_i$ if $u$ is ancestral to several leaves.)
The marked tree $[S, (Y_u)]$ is the $k$-spine tree encoded by the r.v.\
$(W_1,\dots,W_{k-1})$ and $(X_1,\dots,X_k)$.

\paragraph{Construction of $\bar{\QQ}^{k,N}_x$.} 
Let $[S, (Y_u)]$ be the $k$-spine tree constructed above. We attach to
$S$ some subtrees distributed as $\PP_x$ to define a larger marked tree
$[T, (X_u)]$. This yields a random tree with $k$ spines originated from
the $k$ leaves of $S$ at generation $N$. The distribution of these random
variables will be denoted by $\bar{\QQ}^{k,N}_x$.

To construct $T$ from the spine, we first specify the number of subtrees
that need to be attached to each vertex $u$ of the spine. We will
distinguish the degree of a vertex in $S$ and that in the larger tree
$T$. We denote by $d_u$ the number of children in $S$ of $u$. (The degree
of $u$ in $T$ will be denoted by $k_u$ as previously.) We work
conditional on $[S, (Y_u)]$ and assume them to be fixed. Let $(K_u;\, u
\in S, \abs{u} < N)$ be independent variables such that $K_u$ has the
distribution of $K(Y_u)$, biased by its $d_u$-th factorial moment. That
is,
\begin{equation}\label{eq:birth}
    \P(K_u = k) = \frac{k^{(d_u)}}{m_{d_u}(Y_u)} \pi(Y_u, k).
\end{equation}
Among the $K_u$ children of $u$ in $T$, $d_u$ are distinguished as they
correspond to the children of $u$ in $S$. Let $C_{u1} < \dots < C_{ud_u}$
be the labels of these distinguished children, and let us assume that
they are uniformly chosen among the $\binom{K_u}{d_u}$ possibilities. We
can now define the subtree corresponding to $S$ in the larger tree $T$ by
an inductive relabelling of the nodes. For $u \in S$, define $\Psi(u)$
inductively as follows
\[
    \Psi(\emptyset) = \emptyset,\quad \Psi(ui) = \Psi(u)C_{ui},
\]
with corresponding marks
\[
    X_{\Psi(u)} = Y_u.
\]

Finally, let us attach the subtrees to $[S, (Y_u)]$. For $u \in S$,
consider a sequence $[T_{ui}, (X_{ui,v};\, v);\, i \ge 1]$ of i.i.d.\
marked trees with the original distribution $\PP_{X_{ui,\emptyset}}$,
but with random initial mark $X_{ui, \emptyset}$ distributed as
$p(Y_u, \cdot)$. The final tree $T$ is defined as 
\[
    T = \bigcup_{\substack{u \in S\\ \abs{u} < N}}
    \bigcup_{\substack{i \in [K_u]\\ ui \notin \Psi(S)}} 
    \Psi(u)iT_{ui},
\]
and for $v \in T_{ui}$, the mark of $\Psi(u)iv$ is 
\[
    X_{\Psi(u)iv} = X_{ui,v}.
\]
Informally, for each of the $K_u - d_u$ children of $u$ that are not in
$S$, we realize one step of the Markov chain with kernel $(p(x,y);\, x,y
\in E)$ and then attach a whole subtree $T_{ui}$ to that child.

The resulting tree $T$ has $k$ distinguished leaves, 
$\Psi(\ell_1), \dots, \Psi(\ell_k)$, corresponding to the $k$ leaves of
$S$. Let us finally define
\[
    \forall i \le k,\quad  V_i = \Psi(\ell_{\sigma_i})
\]
for an independent uniform permutation $\sigma$ of $\{1,\dots, k\}$. The
distribution of the triple $[T_{[N]}, (X_u), (V_i)]$ is denoted by
$\bar{\QQ}^{k,N}_x$, where $T_{[N]}$ is the restriction of $T$ to the
first $N$ generations.

\subsection{The spinal decomposition theorem}
\label{SS:spinalDecomposition}

Our final objective in this section is to connect $\PP_x$ and
$\bar{\QQ}^{k,N}_x$ to derive our many-to-few formula. We assume $\nu_n >
0$ for all $n \in \{ 0, \dots, N-1 \}$. Recall the expression of
$\Delta_k$ from \eqref{eq:bias}. Our spinal decomposition theorem states
that, if $\PP_x$ is biased by the $k$-th factorial moment of its size at
generation $N$ and $k$ uniformly chosen individuals $(V_1,\dots,V_k)$ are
distinguished from that generation, the corresponding marked tree $[T_N,
(X_u), (V_i)]$ is distributed as $\bar{\QQ}^{k,N}_x$ biased by
$\Delta_k$.

\begin{theorem}[Spinal decomposition] \label{thm:spinal}
    Consider a tree $\tau$ with height $N$ and $k$ distinct vertices
    $(v_1, \dots, v_k) \in \tau_N$. Let $h$ be a harmonic function for
    the branching process with law $(\PP_x;\, x \in E)$. Then, for any
    test function $\phi$, we have
    \[ 
        \PP^N_x\Big[ \indic_{T=\tau}  \phi( X_u;\, u \in \tau )
            \Big]
            = h(x) N^{k-1} k!\,
        \bar{\QQ}^{k,N}_x\Big[ \Delta_k \cdot \indic_{T = \tau} \indic_{V_1=v_1,\dots, V_k=v_k} 
             \phi(X_u;\, u \in \tau) \Big].
    \]
\end{theorem}

\begin{proof}
    It is enough to prove the result for the uniform CPP by
    noting that $\prod_{u\in S}\frac{1}{N\nu_u}$ is the Radon--Nykodim
    derivative of the uniform CPP with respect to the $\nu$-CPP.

    The natural state space for $\PP^N_x$ is the space of all marked
    trees with height at most $N$, that is,
    \[
        \bigcup_{\tau \in \Omega_N} \{ \tau \} \times E^\tau.
    \]
    Using that the offspring distribution on $E$ has a density w.r.t.\
    some measure $\Lambda$, it is clear that $\PP^N_x$ has a density
    w.r.t.\ a dominating measure defined as 
    \[
        \Lambda^N\Big(\bigcup_{\tau \in \Omega_N} \{\tau \} \times \prod_{u \in
        \tau} B_u\Big) = \sum_{\tau \in \Omega_N} \prod_{u \in \tau} \Lambda(B_u)
    \]
    which is given by
    \begin{equation} \label{eq:density}
        P_x^N(\tau, (x_u)) = \prod_{\substack{u \in \tau\\ \abs{u} < n}} 
        \pi(x_u, k_u) \prod_{i=1}^{k_u} p(x_u, x_{ui})
    \end{equation}
    where $k_u$ stands for the number of children of $u$. 

    Let $\spine$ denote the subtree spanned by
    $(v_1,\dots,v_k)$, that is
    \[
        \spine = \bigcup_{i = 1}^k \{ w : w \preceq v_i\}.
    \]
    We can decompose \eqref{eq:density} into a product on $\spine$ and on
    the subtrees attached to $\spine$. The branching property shows that
    \begin{align*}
        P_x^N(\tau, (x_u)) &= \prod_{\substack{u \in \spine\\ \abs{u} <
        N}} \pi(x_u, k_u)
        \bigg[\prod_{\substack{i \in [k_u]\\ ui \notin \spine}} p(x_u, x_{ui})
        P_{x_{ui}}\big(\tau_{ui}, (x_{uiv})_v\big) 
        \times \prod_{\substack{i \in [k_u]\\ ui \in \spine}} p(x_u, x_{ui}) \bigg] \\
        &= \bigg[ \prod_{\substack{u \in \spine\\ \abs{u} < N}} \pi(x_u, k_u)
            \prod_{\substack{i \in [k_u]\\ ui \notin \spine}} p(x_u, x_{ui})
        P_{x_{ui}}(\tau_{ui}, (x_{uiv})_v)  \bigg]
        \prod_{u \in \spine \setminus \{\emptyset\}} p(x_{\an{u}}, x_u).
    \end{align*}

    For $u \in \spine$, let $d_u$ denote the number of children of $u$
    that belong to $\spine$, that is,
    \[
        d_u = \Card \{ i : ui \in \spine \}.
    \]
    Let us make the following change in the
    previous equality
    \[
        \pi(x_u, k_u) = \frac{1}{k_u^{(d_u)}} \cdot
        \frac{k_u^{(d_u)} \pi(x_u, k_u)}{m_{d_u}(x_u)} \cdot m_{d_u}(x_u).
    \]
    Let us also write the second term in the product as
    \begin{align*}
        \prod_{u \in \spine \setminus \{\emptyset\}}
        p(x_{\an{u}}, x_u)
        &=
        \prod_{u \in \spine \setminus \{\emptyset\}}
        \frac{h(x_u) m(x_{\an{u}}) p(x_{\an{u}}, x_u)}{h(x_{\an{u}})} 
        \cdot \frac{h(x_{\an{u}})}{m(x_{\an{u}})h(x_u)} \\
        &=
        \prod_{u \in \spine \setminus \{\emptyset\}}
        q(x_{\an{u}}, x_u)
        \cdot
        \prod_{\substack{u \in \spine\\ \abs{u} < N}} \frac{h(x_u)^{d_u-1}}{m(x_u)^{d_u}}
        \cdot
        \prod_{i=1}^k \frac{1}{h(x_{v_i})}
        \cdot h(x_\emptyset).
    \end{align*}

    Putting both expressions together, we obtain that 
    \begin{align*}
        \frac{1}{h(x_{\emptyset})} 
        &P_{x}(\tau, (x_u)) 
        = 
        \prod_{i = 1}^k \frac{1}{h(x_{v_i})} 
        \prod_{\substack{u \in \spine\\ d_u > 1}}
        h(x_u)^{d_u - 1} \frac{m_{d_u}(x_u)}{d_u!\,m(x_u)^{d_u}}
        \times 
        \prod_{u \in \spine \setminus \{\emptyset\}}
        q(x_{\an{u}}, x_u) \\
        &\times \prod_{\substack{u \in \spine\\ \abs{u} < N}} \bigg[ \frac{d_u!}{k_u^{(d_u)}} \cdot
        \frac{k_u^{(d_u)} \pi(x_u, k_u)}{m_{d_u}(x_u)} 
        \cdot \prod_{ui \notin \spine} p(x_u, x_{ui}) P_{x_{ui}}\big(\tau_{ui}, (x_{uiv})_v\big) \bigg].
    \end{align*}
    The result now follows upon identifying each term in this product.
    The first term is $\Delta_k$. 
    The second term is the density of the marks $(x_u;\, u \in
    \spine)$ along the $k$-spine. The last product is made of three terms.
    The first is the probability that the $d_u$ children of $u$ that
    belong to $\spine$ have a given birth rank. The second is the
    probability that the final degree of $u$ is $k_u$ given that it has
    $d_u$ children in $\spine$ (see \eqref{eq:birth}). The last is the density of the marked
    trees attached to $u$. The $N^{k-1}k!$ term in the statement of
    theorem is simply the probability of observing a given ultrametric
    tree and labeling of the leaves. 
\end{proof}

\begin{proof}[Proof of Proposition~\ref{prop:manyToFew}]
    Let $\tau$ be some fixed tree and $v_1,\dots,v_k$ be distinct
    vertices at height $N$ of $\tau$. Using Theorem~\ref{thm:spinal} 
    yields
    \begin{multline*}
        \PP^N_x\Big[ \indic_{T = \tau}
            \phi\big( d_T(v_i,v_j), X_{v_i};\, i,j \le k \big)
        \Big]  \\
        = 
        h(x) N^{k-1} k!\, \bar{\QQ}^{k,N}_x\Big[ 
            \Delta_k \cdot \indic_{T = \tau}
            \indic_{V_1=v_1,\dots,V_k=v_k}
            \phi\big( d_T(V_i, V_j), Y_{V_i};\, i,j \le k\big)
        \Big].
    \end{multline*}
    Summing over all $(v_1,\dots,v_k)$ first, then over all $\tau$, and
    recalling that $V_i = \ell_{\sigma_i}$ for an independent uniform
    permutation $\sigma$ of $\{1,\dots,k\}$ proves the result.
\end{proof}

\section{Convergence of marked branching processes}
\label{S:GromovWeak}

\subsection{The marked Gromov-weak topology}
\label{SS:gromovWeak}

Deriving the scaling limit of the genealogy and types in a branching
process requires one to envision it as a random marked metric measure space.
In this work we equip the set of all such spaces with the marked
Gromov-weak topology \cite{depperschmidt_marked_2011}. This section is a
brief remainder of the basic properties of this topology, a more thorough 
account can be found in \cite{depperschmidt_marked_2011,
greven_convergence_2006}. We do not restrict our attention to trees and
try to follow as much as possible the notation in \cite{depperschmidt_marked_2011}, 
so that some notation in this section might be inconsistent with the rest
of the paper.

Let $(E, d_E)$ be a fixed complete separable metric space, referred to as
the \emph{mark space}. In our application, $E = [0, \infty)$ is
endowed with the usual distance on the real line. A \emph{marked metric
measure space} (mmm-space for short) is a triple $[X, d, \mu]$, where
$(X,d)$ is a complete separable metric space, and $\mu$ is a finite
measure on $X \times E$.

To define a topology on the set of mmm-spaces, for each $k \ge 1$
consider the map
\[
    R_k \colon 
    \begin{cases}
        (X \times E)^k \to \R_+^{k^2} \times E^k \\
        \big( (x_i, u_i);\, i \le k \big) \mapsto 
        \big( d(x_i, x_j), u_i ;\, i,j \le k \big)
    \end{cases}
\]
that maps $k$ points in $X \times E$ to the matrix of pairwise distances and
vector of marks. We denote by $\nu_{k, X} = \mu^{\otimes k} \circ R_k^{-1}$, the
$k$-th \emph{marked distance matrix distribution} of $[X, d, \mu]$, which is the
pushforward of $\mu^{\otimes k}$ by the map $R_k$. (Note that $\mu$ is
not necessarily a probability distribution.) For some $k \ge 1$ and some
continuous bounded test function 
\[
    \phi \colon \R_+^{k^2} \times E^k \to \R
\]
let us define a functional
\begin{equation} \label{eq:polynomials}
    \Phi\big( X, d, \mu \big) = \Angle{\nu_{k,X}, \phi}.
\end{equation}
Functionals of the previous form are called \emph{polynomials}
($k$ is the \emph{degree} or \emph{order} of the polynomial), and
the set of all polynomials, obtained by varying $k$ and $\phi$, is denoted by $\Pi$.

\begin{definition}
    The marked Gromov-weak topology is the topology on mmm-spaces induced
    by $\Pi$. A random mmm-space is a r.v.\ with values in the set of
    (equivalence classes of) mmm-spaces, endowed with the Gromov-weak
    topology and the associated Borel $\sigma$-field.
\end{definition}

\begin{remark}
    Formally, the marked Gromov-weak topology should be defined on
    equivalence classes of mmm-spaces, where two spaces belong to the
    same class \tiff there is a measure preserving isometry between the
    supports of their measures that also preserves marks, see
    \cite[Definition~2.1]{depperschmidt_marked_2011}. This distinction
    has little consequences in practice so that we often omit it.
\end{remark}

There is a unique equivalence class of all mmm-spaces with a null
sampling measure, which acts as the null mmm-space and that we denote
by $\mathbf{0}$. It follows from the definition of the Gromov-weak
topology that a sequence of mmm-spaces $([X_n, d_n, \mu_n];\, n \ge
1)$ converges to $\mathbf{0}$ \tiff $\mu_n(X_n \times E) \to 0$.
If $[X,d,\mu]$ is a \emph{random} mmm-space, the expectation of a
polynomial evaluated at $[X,d,\mu]$, namely $\E[\Phi(X,d,\mu)]$, is
called a \emph{moment} of $[X,d,\mu]$.

\begin{remark}[Polar decomposition]
    An mmm-space $[X,d,\mu] \ne \mathbf{0}$ can be seen as a pair $(\bar{\mu},
    [X,d,\hat{\mu}])$ where $\bar{\mu} = \mu(X \times E)$ is the total mass of
    $\mu$ and $\hat{\mu} = \mu / \bar{\mu}$ is the renormalized
    probability measure. This is the so-called \emph{polar decomposition}
    of $[X,d,\mu]$ \cite{depperschmidt2019treevalued}. The space of all
    polar decompositions is naturally endowed with the product topology,
    where the space of all \emph{probability} mmm-spaces is endowed with
    the more standard marked Gromov-weak topology restricted to
    probability mmm-spaces \cite{depperschmidt_marked_2011}. It is not
    hard to see that the map taking non-null mmm-spaces to their polar
    decompositions is an homeomorphism.

    An important consequence of this remark is that the convergence in
    distribution of a sequence of mmm-spaces $[X_n, d_n, \mu_n]$ implies
    that of $[X_n, d_n, \hat{\mu}_n]$, provided that the limit mmm-space is
    a.s.\ non-null. In particular, for ultrametric spaces, it implies the
    convergence in distribution of the genealogy of $k$ individuals
    sampled from $[X_n ,d_n ,\mu_n]$ according to $\hat{\mu}_n$.
\end{remark}

Many properties of the marked Gromov-weak topology are derived in
\cite{depperschmidt_marked_2011} under the further assumption that 
$\mu$ is a probability measure. Relaxing this assumption to account for
finite measures is quite straightforward but requires some caution, as
the total mass of $\mu$ can now drift to zero or infinity. In particular,
the following result shows that $\Pi$ forms a convergence determining
class only when the limit satisfies a moment condition, which is a
well-known criterion for a real variable to be identified by its moments,
see for instance \cite[Theorem~3.3.25]{durrett_probability_2019}. 
This result was already stated for metric measure spaces without marks
in \cite[Lemma~2.7]{depperschmidt2019treevalued}. 

\begin{proposition} \label{lem:convDetermining}
    Suppose that $[X,d,\mu]$ is a random mmm-space verifying 
    \begin{equation} \label{eq:momentCondition}
        \limsup_{p \to \infty} \frac{\E[\mu(X\times E)^p]^{1/p}}{p} < \infty.
    \end{equation}
    Then, for a sequence $[X_n, d_n, \mu_n]$ of random mmm-spaces to
    converge in distribution for the marked Gromov-weak topology to
    $[X,d,\mu]$ it is sufficient that
    \[
        \lim_{n \to \infty} \E\big[ \Phi\big(X_n, d_n, \mu_n\big) \big]
        = \E\big[ \Phi\big(X, d, \mu\big) \big]
    \]
    for all $\Phi \in \Pi$.
\end{proposition}

\begin{proof}
    Let us prove this result carefully. Fix a polynomial $\Phi$ of degree
    $k$ associated to a non-negative continuous bounded functional $\phi$.
    Recall the notation $\bar{\mu}_n$ for the total mass of $[X_n, d_n,
    \mu_n]$, which is a r.v.\ with values in $[0, \infty)$, and
    $\hat{\mu}_n = \mu_n / \bar{\mu}_n$. Introduce a new measure
    $M_n^\Phi$ on $[0, \infty)$ such that for any continuous bounded
    function $f \colon [0, \infty) \to \R$
    \[
        \angle{M_n^\Phi, f} = \E[ f(\bar{\mu}_n) \Phi(X_n, d_n, \mu_n) ].
    \]
    The key observation is now that, applying Fubini's theorem, $[X, d, \mu] \mapsto \bar{\mu}^p
    \Phi(X,d,\mu)$ is again a polynomial of the form
    \eqref{eq:polynomials} (of degree $p+k$). Therefore, our assumption
    entails that, for any integer $p \ge 0$,
    \[
        \int_0^\infty x^p M_n^\Phi(\diff x)
        = \E[ \bar{\mu}_n^p \Phi(X_n, d_n, \mu_n) ]
        \goesto{n \to \infty} 
        \E[ \bar{\mu}^p \Phi(X, d, \mu) ]
        = \int_0^\infty x^p M^\Phi(\diff x)
    \]
    where we have defined $M^\Phi$ is a similar way to $M^\Phi_n$ using
    the limiting random variable $[X,d,\mu]$. Now, the usual method of
    moments on $[0, \infty)$, see for instance
    \cite[Theorem~3.3.26]{durrett_probability_2019}, entails that for any
    continuous bounded function $f$ 
    \begin{equation} \label{eq:convDetermining1}
        \E[ f(\bar{\mu}_n) \Phi(X_n, d_n, \mu_n) ]
        = \int_0^\infty f(x) M_n^\Phi(\diff x)
        \goesto{n \to \infty} 
        \int_0^\infty f(x) M^\Phi(\diff x)
        = \E[ f(\bar{\mu}) \Phi(X, d, \mu) ].
    \end{equation}
    We have used that $M^\Phi$ fulfills the moment growth condition
    of \cite[Theorem~3.3.26]{durrett_probability_2019} since
    \[
        \int_0^\infty x^p M^\Phi(\diff x) \le \norm{\phi}_\infty
        \E[\mu(X\times E)^{p+k}]
    \]
    and \eqref{eq:momentCondition} holds. By taking linear combinations,
    \eqref{eq:convDetermining1} holds for any polynomials, not only
    non-negative ones.

    Let $f \colon [0, \infty) \to \R$ be continuous bounded and have its
    support bounded away from $0$. Since $x \mapsto f(x) / x^k$ is
    continuous bounded, applying \eqref{eq:convDetermining1} to this map
    and using that $\Phi(X_n, d_n, \mu_n) = \bar{\mu}_n^k \Phi(X_n, d_n,
    \hat{\mu}_n)$ shows 
    \[
        \E[ f(\bar{\mu}_n) \Phi(X_n, d_n, \hat{\mu}_n) ]
        \goesto{n \to \infty} 
        \E[ f(\bar{\mu}) \Phi(X, d, \hat{\mu}) ].
    \]
    Standard arguments show that the above convergence also holds for
    $f(x) = g(x) \indic_{\{x \ge \epsilon\}}$ for any continuous bounded
    $g$ and $\epsilon > 0$ such that $\P(\bar{\mu} = \epsilon) = 0$.
    Since \cite[Theorem~5]{depperschmidt_marked_2011} ensures that
    polynomials are convergence determining on mmm-spaces with a
    \emph{probability} sampling measure, we can use
    \cite[Proposition~4.6, Chapter 3]{ethier_1986} to obtain that for any
    continuous bounded functional $F$ on the space of mmm-spaces, 
    \begin{equation} \label{eq:convDetermining2}
        \E[ F(X_n, d_n, \mu_n) \indic_{\{\bar{\mu}_n \ge \epsilon\}} ]
        \goesto{n \to \infty} 
        \E[ F(X, d, \mu) \indic_{\{\bar{\mu} \ge \epsilon\}} ].
    \end{equation}
    (Here we have applied the result to the polar decomposition of the
    mmm-space, and used that the polar decomposition defines an 
    homeomorphism.)

    To end the proof, we note that by Portmanteau's theorem
    \begin{equation} \label{eq:convDetermining3}
        \P(\bar{\mu} < \epsilon) \le 
        \liminf_{n \to \infty} \P(\bar{\mu}_n < \epsilon) 
        \le \limsup_{n \to \infty} \P(\bar{\mu}_n \le \epsilon) 
        \le \P(\bar{\mu} \le \epsilon).
    \end{equation}
    Finally, we write 
    \[
        \E[F(X_n, d_n, \mu_n)] = 
        \E[F(X_n, d_n, \mu_n) \indic_{\{\bar{\mu}_n \ge \epsilon\}} ]
        + \E[F(X_n, d_n, \mu_n) \indic_{\{\bar{\mu}_n < \epsilon\}} ]
    \]
    take a limit $n \to \infty$ first, then $\epsilon \to 0$, and use
    \eqref{eq:convDetermining2} to estimate the first term and 
    \eqref{eq:convDetermining3} to control the second one to obtain
    \[
        \E[F(X_n, d_n, \mu_n)]
        \goesto{n \to \infty}
        \E[F(X, d, \mu)]
    \]
    which is the desired result.
\end{proof}

\paragraph{Bi-metric measure spaces.} The branching process with
recombination is naturally endowed with two metrics: the genealogical
distance and the chromosomic distance. Therefore, for the purpose of this
application only, let us say that $[X, d, D, \mu]$ is a \emph{marked bi-metric
measure space} if both $d$ and $D$ are metric that make $(X, d)$ and $(X,
D)$ Polish spaces, and if $\mu$ is a finite measure on $X \times E$,
where $X$ is endowed with the $\sigma$-field induced by reunion of the
open balls of $d$ and $D$.

A polynomial of a marked bi-metric measure space is a functional of the
form
\begin{equation} \label{eq:polynomialsBimetric}
    \int_{(X \times E)^k} 
    \phi\big( d(x_i, x_j), D(x_i, x_j), u_i;\, i,j \le k\big)
    \mu^{\otimes k}(\diff (x_1, u_1) ,\dots, \diff (x_k, u_k))
\end{equation}
for some $k$ and some $\phi$. Accordingly we define the Gromov-weak
topology for these spaces as the topology induced by the polynomials.
It is straightforward to check that all the results stated for mmm-spaces  
carry on to marked bi-metric measure spaces, up to replacing the
polynomials in \eqref{eq:polynomials} by that in
\eqref{eq:polynomialsBimetric}.

\subsection{Convergence of ultrametric spaces}

Using Proposition~\ref{lem:convDetermining} requires one to have prior knowledge of
the limit $[X,d,\mu]$. A stronger version of this result would be that
the convergence of each $(\E[ \Phi(X_n,d_n,\mu_n) ];\, n \ge 1)$ implies
the existence of a random mmm-space to which $(X_n, d_n, \mu_n)$
converges in distribution (under a moment condition similar to
\eqref{eq:momentCondition}). Such a result cannot hold in the current
formulation of the marked Gromov-weak topology. This is a consequence of
the fact that some limits of distance matrix distributions cannot be
expressed as the distance matrix distribution of a \emph{separable}
metric space, see for instance \cite[Example~2.12 (ii)]{greven_convergence_2006}.
To overcome this issue, it is necessary to relax the separability
assumption in the definition of an mmm-space.

Deriving a meaningful extension of the Gromov-weak topology to
non-separable metric spaces is not a straightforward task, since it
raises many measure theoretic difficulties. However, when restricting our
attention to genealogies, as is the purpose of this work, the specific
tree structure of these objects can be used to define such an extension.
We follow the framework introduced in \cite[Section~4]{foutelrodier2019exchangeable}, 
but see also \cite{gufler_representation_2016}. The results contained in
this section are not necessary for the analysis of the branching process
with recombination and can be possibly skipped.

\begin{definition}[Marked UMS, \cite{foutelrodier2019exchangeable}]
    \label{def:UMS}
    A marked ultrametric measure space (marked UMS) is a collection 
    $[U, d, \mathscr{U}, \mu]$ where $\mathscr{U}$ is a $\sigma$-field and
    \begin{enumerate}[(i)]
        \item \label{it:UMS1} The metric $d$ is $\mathscr{U} \otimes
            \mathscr{U}$-measurable and is an ultrametric:
            \[
                \forall x,y,z \in X,\quad d(x,y) \le \max \{ d(x,z), d(z,y) \}.
            \]
        \item \label{it:UMS2} The $\sigma$-field  $\mathscr{U}$ verifies:
            \[
                \sigma\big( B(x,t);\, x \in U, t > 0\big) \subseteq \mathscr{U} \subseteq \mathscr{B}(U)
            \]
            where $B(x,t)$ is the open ball of radius $t$ and center $x$,
            and $\mathscr{B}(U)$ is the Borel $\sigma$-field associated
            to $(U, d)$;
        \item The measure $\mu$ is a finite measure on $U \times E$,
            defined on the product $\sigma$-field $\mathscr{U} \otimes
            \mathscr{B}(E)$.
    \end{enumerate}
\end{definition}

\begin{remark}
    While this definition might be surprising at first sight, note that
    if $(U,d)$ is separable and ultrametric, points \eqref{it:UMS1} and
    \eqref{it:UMS2} of the definition are fulfilled when $\mathscr{U}$ is
    chosen to be the usual Borel $\sigma$-field. Therefore, a separable
    marked UMS in the sense of Definition~\ref{def:UMS} is an ultrametric
    mmm-space in the sense of Section~\ref{SS:gromovWeak}. When no
    $\sigma$-field is prescribed, $\mathscr{U}$ is assumed to be the Borel
    $\sigma$-field. Using a naive definition of a marked UMS as a
    complete metric space with a finite measure on the corresponding
    Borel $\sigma$-field raises some deep measure theoretic issues
    related to the Banach--Ulam problem, that are avoided by
    Definition~\ref{def:UMS}, see
    \cite[Section~4]{foutelrodier2019exchangeable} for a discussion.
\end{remark}

Point~\eqref{it:UMS1} of the above definition ensures that each map $R_k$
is measurable, so that we can define the marked distance matrix
distribution $\nu_{k,U}$ and the polynomials $\Phi(U, d, \mathscr{U},
\mu)$ of a marked UMS $(U,d,\mathscr{U},\mu)$ as in the previous section.
Analogously to mmm-spaces, we define the marked Gromov-weak topology on
the set of marked UMS as the topology induced by the set of polynomials.

\begin{remark}
    Again, for the topology to be separated we need to work with
    equivalence classes of marked UMS. For non-separable spaces, the
    correct notion of equivalence is that of weak isometry provided
    in \cite[Definition~4.11]{foutelrodier2019exchangeable}. We do not
    make the distinction between marked UMS and their equivalence class
    in practice.
\end{remark}

We can now state a stronger version of Proposition~\ref{lem:convDetermining}
for ultrametric spaces. In the statement of the theorem we will
need a mild tightness conditions. For a marked UMS $[U,d,\mathscr{U},
\mu]$, define the maps $r$ and $\pi_E$ as 
\[
    \forall (x_1,u_1), (x_2,u_2) \in U \times E,\quad
    r( (x_1,u_1), (x_2, u_2)) = d(x_1, x_2),\quad
    \pi_E( (x_1,u_1) ) = u_1
\]
and the corresponding pushforward measures
\[
    w_U = \mu^{\otimes 2} \circ r^{-1}, \qquad m_U = \mu \circ \pi_E^{-1}.
\]
If $[U,d,\mathscr{U},\mu]$ is random, these are random measures, and we
denote their intensity measures by $\E[w_U]$ and $\E[m_U]$, which are
deterministic measures on $\R_+$ and $E$ respectively.

\begin{theorem} \label{thm:convUMS}
    Let $(U_n, d_n, \mathscr{U}_n, \mu_n)$ be a sequence of random marked UMS such that
    for any polynomial $\Phi \in \Pi$,
    \[
        \lim_{n \to \infty} \E\big[ \Phi\big( U_n, d_n, \mathscr{U}_n, \mu_n \big)
        \big]
    \]
    exists, and fulfill (compare with \eqref{eq:momentCondition})
    \begin{equation} \label{eq:momentUMS}
        \limsup_{p \to \infty} \lim_{n \to \infty} 
        \frac{\E[ \mu_n(U_n \times E)^p]^{1/p} }{p} < \infty.
    \end{equation}
    Suppose also that the sequences $(\E[w_{U_n}];\, n \ge 1)$ and
    $(\E[m_{U_n}];\, n \ge 1)$ are relatively compact, as measures on
    $\R_+$ and $E$. Then there exist a random marked UMS, $[U, d,
    \mathscr{U}, \mu]$ such that $(U_n, d_n, \mathscr{U}_n, \mu_n)$
    converges to that limit in the marked Gromov-weak topology. Moreover
    the limit is characterized
    by
    \[
        \E\big[ \Phi\big( U, d, \mathscr{U}, \mu \big)\big]
        = \lim_{n \to \infty}
        \E\big[ \Phi\big( U_n, d_n, \mathscr{U}_n, \mu_n \big)\big].
    \]
\end{theorem}

\begin{remark}
    The previous result suggests the following simple method to
    prove convergence in distribution in the (usual) sense of separable
    ultrametric mmm-spaces. First prove that the conditions of
    Theorem~\ref{thm:convUMS} are fulfilled, then check that the limiting
    marked UMS is a.s.\ separable. The two compactness conditions
    on $(\E[w_{U_n}])$ and $(\E[m_{U_n}])$ ensure, in combination with
    the convergence of the moments, that the sequence of mmm-spaces is
    tight. Compare this to checking, on top of the previous assumptions,
    the tightness criterion in \cite[Theorem~2 (ii)]{greven_convergence_2006} 
    that ensures that no mass of the sampling measure is accumulating on
    isolated points. This condition is not needed here because we have
    enlarged the state space of mmm-spaces to include non-separable
    metric spaces.
\end{remark}

The proof of the above result is based on a characterization of all
exchangeable ultrametric matrices. We call a random pair $(d_{ij};\, i,j
\ge 1)$ and $(Y_i;\, i \ge 1)$ a \emph{marked exchangeable ultrametric
matrix} if 
\begin{itemize}
    \item each $Y_i$ has values in $E$;
    \item $(d_{ij};\, i, j \ge 1)$ is a.s.\ an ultrametric on $\N$;
    \item its distribution is invariant by the action of any permutation
    $\sigma$ of $\N$ with finite support:
    \[
        \big( (d_{\sigma_i\sigma_j};\, i, j \ge 1), 
            (Y_{\sigma_i};\, i \ge 1) \big)
        \overset{\mathrm{(d)}}{=}
        \big( (d_{ij};\, i, j \ge 1), 
        (Y_i;\, i \ge 1) \big).
    \]
\end{itemize}

A typical way to obtain such an ultrametric matrix is to consider an
i.i.d.\ sample $(X_i, Y_i;\, i \ge 1)$ from a marked UMS
$(U,d,\mathscr{U},\mu)$ with $\mu(U) = 1$ a.s., and define 
\begin{equation} \label{eq:sampleUMS}
    \forall i, j \ge 1,\quad d_{ij} = d(X_i, X_j).
\end{equation}
The next result shows that all exchangeable marked matrices are
obtained in this way. It can be seen as a version of Kingman's
representation theorem of exchangeable partitions
\cite{kingman_representation_1978} for ultrametric matrices. 

\begin{theorem}[\cite{foutelrodier2019exchangeable}]
    \label{thm:deFinettiUMS}
    Let $(d_{ij};\, i,j \ge 1)$ and $(Y_i;\, i \ge 1)$ be an exchangeable
    marked ultrametric matrix. There exists a random marked probability
    UMS $[U, d, \mathscr{U}, \mu]$ (that is, $\mu(U)=1$ a.s.)
    such that the exchangeable marked ultrametric matrix obtained by
    sampling from it as in \eqref{eq:sampleUMS} is distributed as
    $(d_{ij};\, i,j \ge 1)$ and $(Y_i;\, i \ge 1)$. Moreover this marked
    UMS is unique in distribution.
\end{theorem}

\begin{proof}
    This result is a straightforward extension of
    \cite[Theorem~1.8]{foutelrodier2019exchangeable} that deals with the
    case without marks. To guide the reader, let us mention the crucial
    modification that need to be made. The proof relies on encoding some
    marginals of $(d_{ij};\, i,j \ge 1)$ as an exchangeable sequence of 
    r.v.\ $(\xi^{(0)}_i, \dots, \xi_i^{(p)};\, i \ge 1)$ in $[0, 1]^p$ and using a
    de~Finetti-type argument, see
    \cite[Appendix~B]{foutelrodier2019exchangeable}. The same argument
    should be applied to the exchangeable sequence of r.v.\ 
    $(\xi^{(0)}_i,\dots,\xi^{(p)}_i, Y_i;\, i \ge 1, )$.
\end{proof}

\begin{proof}[Proof of Theorem~\ref{thm:convUMS}]
    We prove the result by a tightness and uniqueness argument. 
    To prove tightness, we embed the space of marked UMS into a space of
    measures, using the marked distance matrices, and use known tightness
    arguments for random measures. More precisely, the map $\iota \colon
    [U, d, \mathscr{U}, \mu] \mapsto (\nu_{k,U};\, k \ge 1)$ is an
    injection. This is a consequence of the uniqueness part of
    Theorem~\ref{thm:deFinettiUMS}. For each $k \ge 1$, $\nu_{k,U}$ lives
    in the space of finite measures on $\R_+^{k^2} \times E^k$, which can
    be endowed with the weak topology. If the space of sequences
    $(\nu_{k,U};\, k \ge 1)$ is endowed with the product topology, it
    follows readily from the definition of the Gromov-weak topology that
    $\iota$ is a homeomorphism from the space of marked UMS to its image.
    We claim that the image of $\iota$ is closed in this product
    topology. If this is the case, the space of marked UMS is
    homeomorphic to a closed subset of the space of sequences of
    measures, and clearly, 
    \[
        \text{$(U_n, d_n, \mathscr{U}_n, \mu_n;\, n \ge 1)$ is tight}
        \iff 
        \forall k \ge 1,\;
        \text{$(\nu_{k,U_n};\, n \ge 1)$ is tight},
    \]
    where in the right-hand side each $\nu_{k,Un}$ is a random measure,
    and tightness is with respect to the weak topology. For the collection 
    of random measures $(\nu_{k,U_n};\, n \ge 1)$ to be tight it is
    sufficient that the collection of intensity measures
    $(\E[\nu_{k,U_n}];\, n \ge 1)$ is relatively compact
    \cite[Lemma~3.2.8]{Dawson93}. For $y = ((d_{ij}), u_i;\, i,j \le k)
    \in \R^{k^2} \times E^k$, let $p_{ij}(y) = d_{ij}$ and $p'_i(y) =
    u_i$ be the projection maps. It is sufficient to show that the
    pushforward of $\E[\nu_{k,U_n}]$ through each projection is
    relatively compact. By definition, for a Borel set $A \subseteq \R_+$
    and $i \ne j$, by exchangeability,
    \[
        \big(\E[\nu_{k,U_n}] \circ p_{ij}^{-1}\big)(A) = 
        \E\Big[ \int_{(U_n \times E)^k} \indic_{\{d_n(x_i,x_j) \in A\}}
            \,\diff \mu_n^{\otimes k}( (x_i, u_i)_i )
        \Big]
        = \E\big[ \bar{\mu}_n^{k-2} w_{U_n}(A) \big].
    \]
    The relative compactness of $(\E[\nu_{k,U_n}] \circ p_{ij}^{-1};\, n
    \ge 1)$ now follows from that of $(\E[w_{U_n}];\, n \ge 1)$ and from
    the uniform integrability of $\bar{\mu}_n^{k-2}$. In a similar way,
    for a Borel set $B \in E$, 
    \[
        \big(\E[\nu_{k,U_n}] \circ {p'_i}^{-1}\big)(B) = 
        \E\Big[ \int_{(U_n \times E)^k} \indic_{\{u_i \in B\}}
            \,\diff \mu_n^{\otimes k}( (x_i, u_i)_i )
        \Big]
        = \E\big[ \bar{\mu}_n^{k-1} m_{U_n}(A) \big].
    \]
    The desired compactness follows form that of $(\E[m_{U_n}];\, n \ge
    1)$ and from the uniform integrability of $\bar{\mu}_n^{k-1}$.
    
    We now go back to our claim that the image of $\iota$ is closed. For
    each $k, n \ge 1$, let $\nu_{k,n}$ be the $k$-th marked distance
    matrix distribution of some marked UMS $[U_n, d_n, \mathscr{U}_n, \mu_n]$, and
    assume that it converges as $n \to \infty$ to some $\nu_k$. We need to
    show that the limiting sequence of distance matrices can be obtained
    by sampling from a marked UMS. We can assume without loss of
    generality that $\nu_k \neq 0$. Let $\hat{\nu}_k$ be the probability
    measure obtained by renormalizing $\nu_k$, and define similarly
    $\hat{\nu}_{k,n}$. Since the projection of $\hat{\nu}_{k+1,n}$ on
    $\R_+^{k^2} \times E^k$ is equal to $\hat{\nu}_{k,n}$, the same
    property holds for $\hat{\nu}_{k+1}$ and $\hat{\nu}_k$. Using
    Kolmogorov's extension theorem, we can extend consistently the
    measures $(\hat{\nu}_k;\, k \ge 1)$ to a measure $\hat{\nu}_\infty$
    on $\R_+^{\N \times \N} \times E^{\N}$ whose projections on
    finite-dimensional spaces are given by the measures $(\hat{\nu}_k;\,
    k \ge 1)$. Quite clearly, $\hat{\nu}_\infty$ is the law of a marked
    exchangeable ultrametric matrix. (Exchangeability and almost sure
    ultrametricity hold for a fixed $n$, and pass to the limit.)
    Theorem~\ref{thm:deFinettiUMS} shows that we can find a marked UMS
    $[U, d, \mathscr{U}, \hat{\mu}]$ whose $k$-th marked distance matrix
    distribution is $\hat{\nu}_k$. Denote by $\bar{\mu}$ the limit of the
    total mass of $\nu_{1,n}$, and by $\mu = \bar{\mu} \hat{\mu}$. The
    $k$-th marked distance matrix distribution of the marked UMS $[U, d,
    \mathscr{U}, \mu]$ is $\bar{\mu}^k \hat{\nu}_k = \nu_k$. This proves
    the claim. 

    Finally, we prove uniqueness. Let $[U, d, \mathscr{U}, \mu]$ and
    $[U', d', \mathscr{U}', \mu']$ be two random marked UMS, that are
    limits in distribution of a subsequence of $([U_n, d_n,
    \mathscr{U}_n, \mu_n];\, n \ge 1)$. We want to show that they have
    the same distribution. For any polynomial $\Phi \in \Pi$, since
    $\Phi$ is continuous and $(\Phi(U_n, d_n, \mathscr{U}_n, \mu_n);\, n
    \ge 1)$ is uniformly integrable (it has uniformly bounded moments of
    all orders), the moments of the two limiting marked UMS coincide and
    verify \eqref{eq:momentUMS}, namely, 
    \[
        \lim_{n \to \infty} 
        \E\big[ \Phi\big( U_n, d_n, \mathscr{U}_n, \mu_n\big) \big]
        = \E\big[ \Phi\big( U, d, \mathscr{U}, \mu\big) \big]
        = \E\big[ \Phi\big( U', d', \mathscr{U}', \mu'\big) \big].
    \]
    Introducing the same measure $M^\Phi$ as in the proof of
    Proposition~\ref{lem:convDetermining}, the method of moments on
    $\R_+$ shows that, for any continuous bounded $f \colon \R_+ \to \R$
    and any polynomial $\Phi$,
    \begin{equation*}
        \E\big[ f(\bar{\mu})  \Phi\big(U, d, \mathscr{U}, \mu\big) \big]
        =
        \E\big[ f(\bar{\mu}') \Phi\big(U', d', \mathscr{U}', \mu'\big) \big]
    \end{equation*}
    so that if $f(0) = 0$, we have
    \begin{equation} \label{eq:convUMS1}
        \E\big[ f(\bar{\mu})  \Phi\big(U, d, \mathscr{U}, \hat{\mu}\big) \big]
        =
        \E\big[ f(\bar{\mu}') \Phi\big(U', d', \mathscr{U}', \hat{\mu}'\big) \big].
    \end{equation}
    On the event $\{ \bar{\mu} > 0 \}$, let $(d_{ij}, Y_i;\, i,j \ge 1)$
    be the marked exchangeable ultrametric matrix obtained from an
    i.i.d.\ sample from $[U, d, \mathscr{U}, \hat{\mu}]$, and define
    $(d'_{ij}, Y'_i;\, i,j \ge 1)$ similarly from $[U', d', \mathscr{U}',
    \hat{\mu}']$. The identity \eqref{eq:convUMS1} can be written as 
    \begin{equation*}
        \E\big[ f(\bar{\mu}) \phi\big( d_{ij}, Y_i;\, i,j \le k \big) \big]
        =
        \E\big[ f(\bar{\mu}') \phi\big( d'_{ij}, Y'_i;\, i,j \le k \big) \big].
    \end{equation*}
    This equation shows that $\bar{\mu}$ and $\bar{\mu}'$ have the same
    distribution, and for $\bar{\mu}$-a.e.\ $x$, the law of $(d_{ij}, Y_i;\,
    i,j \le k)$ conditional on $\bar{\mu} = x$ is the same as that of 
    $(d'_{ij}, Y'_i;\, i,j \le k)$ conditional on $\bar{\mu}' = x$.
    The uniqueness part of Theorem~\ref{thm:deFinettiUMS} shows that the
    law of $[U, d, \mathscr{U}, \hat{\mu}]$ conditional on $\bar{\mu} =
    x$ is the same as that of  $[U', d', \mathscr{U}', \hat{\mu}']$
    conditional on $\bar{\mu}'$. Combining this with the fact that 
    $\bar{\mu}$ and $\bar{\mu}'$ have the same distribution and that the
    polar decomposition is a homeomorphism proves that the two marked UMS
    have the same distribution.
\end{proof}

\subsection{Moments of some continuous trees}

In this section we compute the moments of some usual random tree
models, namely CPP trees and $\Lambda$-coalescents, to illustrate the
type of expression that can arise for the limiting mmm-space of 
Proposition~\ref{lem:convDetermining}.

\paragraph{Continuous coalescent point processes.}
Coalescent point process trees are a class of continuous random
trees that correspond to the scaling limit of the genealogy of various
branching processes \cite{duchamps18, lambert2010, popovic2004}. Of
particular interest is the Brownian CPP described in
Section~\ref{SS:mainResults} that corresponds to the scaling limit of
critical Galton--Watson processes, and also corresponds to the limit of
the rescaled genealogy of the branching process with recombination.

Consider a Poisson point process $P$ on $[0, \infty) \times (0, \infty)$,
with intensity $\diff t \otimes \nu(\diff x)$. We make the further
assumptions that 
\[
    \forall x > 0,\quad \nu([x, \infty)) < \infty,
    \qquad \nu((0, \infty)) = \infty.
\]
For some $x_0 > 0$, let $Y$ denote the first atom of $P$ whose second
coordinate exceeds $x_0$, that is,
\[
    Y = \inf \{ t \ge 0 : (t,x) \in P,\; x > x_0 \}. 
\]
The CPP tree at height $x_0$ associated to $\nu$ is the random metric
measure space $[(0, Y), d_P, \Leb]$ with
\[
    \forall x \le y,\quad d_P(x,y) = \sup \{ z : (t,z) \in P,\; x \le z \le y \}.
\]

\begin{proposition} \label{prop:CPPpolynomials}
    Let $[(0,Y), d_P, \Leb]$ be the CPP tree at height $x_0$ associated
    to the measure $\nu$. Then for any continuous bounded function $\phi$
    with associated polynomial $\Phi$, we have
    \[
        \E\big[ \Phi\big((0,Y), d_P, \Leb\big) \big] = \frac{k!}{\nu((x_0,\infty))^k}
        \E\big[ \phi\big( H_{\sigma_i,\sigma_j};\, i,j \le k \big)
        \big]
    \]
    where for $i < j$,
    \[
        H_{i,j} = H_{j,i} = \max \{ H_i, \dots, H_{j-1} \},
    \]
    the r.v.\ $(H_1, \dots, H_{k-1})$ are i.i.d.\ with c.d.f.\  
    \[
        \forall x \in [0, x_0],\quad 
        \P( H_1 \le x ) = \frac{\nu((x_0, \infty))}{\nu((x, \infty))},
    \]
    and $\sigma$ is an independent uniform permutation of $\{1, \dots, k\}$.
\end{proposition}

\begin{proof}
    According to \eqref{eq:polynomials}, we need to study the distance of
    $k$ variables sampled uniformly from $[0, Y]$, after having biased
    $([0, Y], d_P, \Leb)$ by the $k$-th moment of its mass.

    Since $Y$ is independent of the restriction of $P$ to $[0,\infty)
    \times [0, x_0]$, the distribution of $([0, Y], d_P, \Leb)$ biased by
    the $k$-th moment of $Y$ is simply that of $([0, Z], d_P, \Leb)$,
    where $Z$ is distributed as $Y$, biased by its $k$-th moment. Let us
    use the notation $\theta \coloneqq \nu((x_0, \infty))$. It is well
    known that $Z$ follows a $\mathrm{Gamma}(\theta, k+1)$ distribution,
    that is, $Z$ has density
    \[
        \theta^{k+1} \frac{x^k}{k!} e^{-\theta x} \diff x.
    \]
    Conditional $Z$, let $(U_1,\dots, U_k)$ be i.i.d.\ uniform variables
    on $[0, Z]$, and denote by $(U^*_1,\dots, U^*_k)$ their order
    statistics. Let us also denote $U^*_0 = 0$ and $U^*_{k+1} = Z$. It is
    standard that 
    \[
        (U^*_1-U^*_0, U^*_2 - U^*_1, \dots, U^*_{k+1} - U^*_k)
    \]
    are independent exponential variables with mean $1/\theta$. Define 
    \[
        \forall i \le k,\quad H_i = d_P(U^*_i, U^*_{i+1}).
    \]
    As the restriction of $P$ to $[0, \infty) \times [0, x_0]$ is
    independent of the vector $(U^*_0, \dots, U^*_{k+1})$,
    $(H_0,\dots,H_k)$ are i.i.d.\ and distributed as
    \[
        \max \{ x : (t,x) \in P,\, t \le Y \}.
    \]
    The following direct computation shows that this has the required
    distribution,
    \begin{align*}
        \P(H_1 \le x) = \P( P([0,Y] \times (x, x_0]) = 0) 
        &= \int_0^\infty \theta e^{-\theta t} \exp\Big(- t \nu((x, x_0]) \Big) \diff t \\
        &= \frac{\theta}{\theta + \nu((x, x_0])}
        = \frac{\nu((x_0, \infty))}{\nu((x, \infty))}.
    \end{align*}
    It is clear from the definition of $d_P$ that for $i < j$,
    \[
        d_P( U^*_i, U^*_j ) = \max \{ H_i, \dots, H_{j-1} \}.
    \]
    Therefore, if $\sigma$ denotes the unique permutation of $[k]$ such
    that $U_i = U^*_{\sigma_i}$, 
    \begin{align*}
        \E\Big[ \int_{[0,Y]^k} \phi\big( d_P(x_i, x_j);\, i,j \le k\big)
        \diff x_1\dots \diff x_k \Big]
        &= \E\big[ Y^k \big]
        \E\big[ \phi\big( d_P(U_i, U_j);\, i,j \le k\big) \big] \\
        &= \frac{k!}{\theta^k}
        \E\big[ \phi\big( d_P(U_i, U_j);\, i,j \le k\big) \big] \\
        &= \frac{k!}{\theta^k}
        \E\big[ \phi\big( H_{\sigma_i, \sigma_j};\, i,j \le k\big) \big].
        \qedhere
    \end{align*}
\end{proof}

In this work, the scaling limit of the genealogy is given by the Brownian
CPP, which is the CPP with height $1$ associated to the measure 
\[
    \nu(\diff x) = \frac{1}{x^2} \diff x.
\]

\begin{corollary} \label{cor:BrownianCPPpolynomials}
    The moments of the Brownian CPP are given by
    \[
        \E\big[ \Phi\big((0,Y), d_P, \Leb\big) \big] = k!\, 
        \E\big[ \phi\big( H_{\sigma_i,\sigma_j};\, i,j \le k \big)
        \big]
    \]
    where for $i < j$,
    \[
        H_{i,j} = H_{j,i} = \max \{ H_i, \dots, H_{j-1} \},
    \]
    the r.v.\ $(H_1, \dots, H_{k-1})$ are i.i.d.\ uniform on $(0,1)$,
    and $\sigma$ is an independent uniform permutation of $\{1, \dots, k\}$.
\end{corollary}

\begin{proof}
    A direct computation shows that 
    \[
        \nu((1, \infty)) = 1,\qquad \frac{\nu((1,\infty))}{\nu((x,\infty))} = x
    \]
    so that the variables $H_i$ in Proposition~\ref{prop:CPPpolynomials}
    are uniform on $[0, 1]$. 
\end{proof}

\paragraph{Metric measure spaces with independent types.}
In our model and in many other settings, the types in the population
become independent of the genealogy in the limit of large population
size. Typically, this situation arises when the time between the
ancestors of two typical individuals in the population is large, so that
the dynamics of the types along the lineages has time to reach some form
of equilibrium and to forget about its starting point (the type of the
ancestor).

For a mmm-space $[X, d, \mu]$, the independence between the types and the
genealogy corresponds to having a product sampling measure of the form 
$\mu = \mu_X \otimes \mu_E$, where $\mu_X$ is a measure on $X$, and
$\mu_E$ a probability measure on the type space $E$. The moments of
such product mmm-spaces are easily expressed in terms of the (unmarked)
metric measure space $[X, d, \mu_X]$.

\begin{proposition} \label{prop:productSpace}
    Let $[X,d,\mu]$ be a random mmm-space with a sampling measure of the
    form $\mu = \mu_X \otimes \mu_E$, where $\mu_E$ is a deterministic
    probability measure on $E$. Then, for any polynomial $\Phi \in \Pi$,
    we have
    \[
        \E\big[ \Phi\big( X, d, \mu \big) \big]
        = \E\Big[
            \int_{X^k} \phi\big( d(x_i, x_j), Y_i;\, i,j \le k \big) \mu_X^{\otimes k}(\diff x_1, \dots, \diff x_k)
        \Big]
    \]
    where $(Y_1, \dots, Y_k)$ are i.i.d., distributed as $\mu_E$, and
    independent of $[X,d,\mu_X]$.
\end{proposition}

\begin{proof}
    By definition of a polynomial and applying Fubini's theorem for
    a.s.\ all realizations of the random measure,
    \begin{align*}
        \E\big[ \Phi\big( X, d, \mu \big) \big]
        &= \E\Big[
            \int_{(X\times E)^k} \phi\big( d(x_i, x_j), u_i;\, i,j \le k \big) 
            (\mu_X \otimes \mu_E)^{\otimes k}(\diff (x_1, u_1), \dots, \diff (x_k, u_k))
        \Big] \\
        &= \E\Big[
            \int_{E^k} \int_{X^k} \phi\big( d(x_i, x_j), u_i;\, i,j \le k \big) 
            \mu_X^{\otimes k}(\diff x_1, \dots, \diff x_k)
            \mu_E^{\otimes k}(\diff u_1, \dots, \diff u_k)
        \Big] \\
        &= \E\Big[
            \int_{X^k} \phi\big( d(x_i, x_j), Y_i;\, i,j \le k \big) \mu_X^{\otimes k}(\diff x_1, \dots, \diff x_k)
        \Big].
    \end{align*}
\end{proof}

\paragraph{$\Lambda$-coalescents.} A $\Lambda$-coalescent is a process
with values in the partitions of $\N$ such that for any $n$, its
restriction to $\{1, \dots, n\}$ is a Markov process with the following
transitions. When the process has $b$ blocks, any $k$ blocks merge at
rate $\lambda_{b,k}$ where
\[
    \lambda_{b,k} = \int_0^1  x^{k-2} (1-x)^{b-k} \Lambda(\diff x)
\] 
for some finite measure $\Lambda$. These processes were introduced in
\cite{pitman_coalescents_1999, sagitov_1999}, and provide the limit
of the genealogy of several celebrated population models with fixed
population size \cite{mohle2001, schweinsberg2003}.

A $\Lambda$-coalescent can be seen as a random ultrametric space on $\N$.
It is possible to take an appropriate completion of this space to 
define an ultrametric $d_\Lambda$ on $(0, 1)$ that encodes the metric
structure of the coalescent, see \cite[Section~4]{greven_convergence_2006} 
for the separable case, and \cite[Section~3]{foutelrodier2019exchangeable} 
for the general case. More precisely, there exists a random ultrametric
$d_\Lambda$ such that if $(V_i;\, i \ge 1)$ is an independent sequence of
i.i.d.\ uniform r.v.\ on $(0,1)$, and $\Pi_t$ is the partition defined
through the equivalence relation
\[
    i \sim_{\Pi_t} j \iff d_\Lambda(V_i, V_j) \le t,
\]
then $(\Pi_t;\, t \ge 0)$ is distributed as a $\Lambda$-coalescent. 
In particular, this leads to the following expression for the moments
of the metric measure space $[(0,1), d_\Lambda, \Leb]$.

\begin{proposition}
    Let $[(0, 1), d_\Lambda, \Leb]$ be a $\Lambda$-coalescent tree.
    Then
    \[
        \E \big[ \Phi\big( (0,1), d_\Lambda, \Leb \big) \big] 
        = \E \big[ \phi\big( d_{ij};\, i,j \le k \big) \big]
    \]
    where
    \[
        \forall i, j \le k,\quad d_{ij} = \inf \{ t \ge 0 : i \sim_{\Pi_t} j \}
    \]
    for a realization $(\Pi_t;\, t \ge 0)$ of a $\Lambda$-coalescent.
\end{proposition}

\subsection{Relating spine convergence to Gromov-weak convergence}

Let $[T, (X_u)]$ be the random marked tree with distribution $\PP_x$
constructed in Section~\ref{SS:notation}, and let $Z_N = \abs{T_N}$
denote the population size at generation $N$. Recall that $T$ can be
endowed with the graph distance $d_T$, and that $T_N$ denotes the $N$-th
generation of the process. The metric $d_T$ restricted to $T_N$ encodes
the genealogy of the population, and has the simple expression
\[
    \forall u,v \in T_N,\quad d_T(u,v) = N - \abs{u \wedge v}.
\]
Define the mark measure on $T_N \times E$ as 
\[
    \mu_N = \sum_{u \in T_N} \delta_{ (u, X_u) }.
\]
The triple $[T_N, d_T, \mu_N]$ is the mmm-space associated to the
branching process $[T, (X_u)]$. The polynomial of degree $k$
corresponding to a functional $\phi$ can be written as
\[
    \Phi\big( T_N, d_T, \mu_N \big) =
    \sum_{(v_1,\dots,v_k) \in T_N} 
        \phi\big( d_T(v_i, v_j), X_{v_i};\, i,j \le k) \big).
\]

The aim of this section is to provide a general convergence criterion for
a rescaling of the sequence of mmm-spaces $[T_N, d_T, \mu_N;\, N \ge 1]$
that only involves computation on the $k$-spine tree. For each $N \ge 1$,
consider a rescaling parameter $\alpha_N$ for the population size,
$\beta_N \colon E \to E$ for the mark space, and $\gamma_N \colon \R_+
\to \R_+$ for the genealogical distances. We assume that $\gamma_N$ is
increasing so that $\gamma_N\circ d_T$ is also an ultrametric, and that
$\alpha_N \to \infty$.

\begin{theorem} \label{thm:branchingConvergence}
    Suppose that for any $k \ge 1$ and any continuous bounded function $\phi$,
    the sequence
    \begin{equation} \label{eq:spinePolynomial}
        \frac{N^{k-1}}{\alpha_N^k \PP^N_x(Z_N > 0)}
        \QQ_x^{k,N}\Big[ \Delta_k \cdot 
            \phi\big( \gamma_N \circ d_T(i, j), \beta_N \circ X_i(N));\,
            i, j \le k \big) \Big]
    \end{equation}
    converges and that the limit fulfills \eqref{eq:momentUMS}.
    Then there exists a random marked UMS $[U,d,\mathscr{U},\mu]$ such
    that conditional on $Z_N > 0$,
    \[
        \lim_{N \to \infty} \Big[T_N, \gamma_N(d_T), \frac{\mu \circ \beta_N^{-1}}{\alpha_N} \Big]
        = [U, d, \mathscr{U}, \mu]
    \]
    holds in distribution for the marked Gromov-weak topology.
\end{theorem}

\begin{proof}
    According to Theorem~\ref{thm:convUMS} it is sufficient to prove that
    the following moments converge,
    \[
        M_N
        \coloneqq  \frac{\PP^N_x(Z_N > 0)^{-1}}{\alpha_N^k}
        \PP^N_x \Big[ \sum_{(u_1,\dots,u_k) \in T_N} 
            \phi\big( \gamma_N \circ d_T(u_i, u_j),
            \beta_N(X_{u_i});\, i,j \le k \big)
        \Big].
    \]
    Let us denote by
    \[
        \widetilde{M}_N
        \coloneqq  \frac{\PP^N_x(Z_N > 0)^{-1}}{\alpha_N^k}
        \PP^N_x \Big[ \sum_{\substack{(u_1,\dots,u_k) \in T_N\\u_1\ne\dots\ne u_k}} 
            \phi\big( \gamma_N \circ d_T(u_i, u_j),
            \beta_N(X_{u_i});\, i,j \le k \big)
        \Big].
    \]
    By the many-to-few formula, Proposition~\ref{prop:manyToFew}, 
    \[
        \widetilde{M}_N
        = 
        \frac{k!\,N^{k-1}}{\alpha_N^k \PP^N(Z_N > 0)}
        \QQ_x^{k,N}\Big[ \Delta_k \cdot 
            \phi\big( \gamma_N \circ d_T(\sigma_i, \sigma_j), \beta_N
                \circ X_{\sigma_i}(N);\,
            i, j \le k \big) \Big]
    \]
    for an independent uniform permutation $\sigma$ of $[k]$. Taking
    $\phi \equiv 1$, the assumption of the result readily implies that 
    \[
        \PP^N_x\big[ Z_N^k \:\big|\: Z_N > 0 \big]
        = O_N( \alpha_N^k ).
    \]
    Therefore, since
    \[
        \abs{M_N - \widetilde{M}_N} = O_N\Big( \PP^N_x\Big[
            \frac{Z_N^{k-1}}{\alpha_N^k} \:\Big|\: Z_N > 0
    \Big]\Big) \goesto{N \to \infty} 0,
    \]
    the convergence of each $M_N$ follows from that of
    $\widetilde{M}_N$ and the result is proved.
\end{proof}

\subsection{Convergence of the \texorpdfstring{$k$}{k}-spine}
\label{SS:Nconvergence}

Theorem~\ref{thm:branchingConvergence} shows that convergence of the
branching process in the Gromov-weak topology can be deduced from the
convergence of some functionals of the $k$-spine tree. We now provide a
general convergence result for the $k$-spine tree that will be used to
compute the limit of \eqref{eq:spinePolynomial} for the branching process
with recombination.

We work under the measure $\QQ^{k,N}_x$ and define 
\[
    W^N_i = \frac{W_i}{N},\qquad 
    \forall t \ge 0,\quad X^N_i(t) = X_i\big( \floor{Nt} \big).
\]
Since our work involves working under various measures, for a sequence
$(P_n;\, n \ge 1)$ of probability measures and a sequence $(Y_n;\, n \ge
1)$ of r.v., we will use the notation
\[
    Y_n \goesto[P_n]{n \to \infty} Y
\]
to mean that the distribution of $Y_n$, under the measure $P_n$,
converges to the distribution of $Y$.

\begin{assumption}[A\ref{as:Feller}] \label{as:Feller}
    \begin{enumerate}[(i)]
        \item There exists a limiting r.v.\ $W$ such that
            \[
                W^N_1 \goesto[\QQ^{1,N}_x]{N \to \infty} W.
            \]
        \item There exists a limiting Feller process $X$ such that, if
            $X^N_1(0) \to X(0)$, 
            \[
                X^N_1 \goesto[\QQ^{1,N}_x]{N \to \infty} X
            \]
            in the Skorohod topology.
    \end{enumerate}
\end{assumption}

There exist equivalent formulations of the second point involving
generators or semigroups, see for instance 
\cite[Theorem~19.28]{kallenberg_foundations_2002}. In the next result, we
use the notation
\[
    \cat{f; t; g} \colon t \mapsto 
    \begin{cases}
        f(s) &\text{if $s < t$} \\
        g(s-t) &\text{if $s \ge t$}
    \end{cases}
\]
for the concatenation of $f$ and $g$ at time $t$.

\begin{proposition} \label{prop:Nconvergence}
    Suppose that (A\ref{as:Feller}) holds. Then
    \[
        \big( (W^N_1, \dots, W^N_{k-1}), (X^N_1,\dots,X^N_k) \big)
        \goesto[\QQ^{k,N}_x]{n \to \infty}
        \big( (W_1, \dots, W_{k-1}), (X_1,\dots,X_k) \big).
    \]
    where,
    \begin{itemize}
        \item the r.v.\ $(W_1, \dots, W_{k-1})$ are i.i.d.\ copies of the
            limiting r.v.\ $W$;         
        \item $X_1$ is distributed as $X$ started from $x$ and is
            independent of $(W_1,\dots, W_{k-1})$;
        \item for each $i$, conditional on $(W_1, \dots, W_{k-1})$ and
        $(X_1, \dots, X_i)$, 
        \[
            X_{i+1} = \cat{X_i; W_i; X'}
        \]
        where $(X'(t);\, t \ge 0)$ is distributed as $X$ started from 
        $X_i(W_i)$.
    \end{itemize}
\end{proposition}

\begin{proof}
    Let us work inductively, and assume that the convergence holds for some
    $k \ge 1$. Let  $\widetilde{X}^N_{k+1}$ be distributed as $X$, started
    from $X^N_k(W^N_k)$.

    Obviously, $W^N_k$ converges to $W_k$, a copy of $W$ independent of
    $(W_1,\dots,W_{k-1})$ and of $(X_1,\dots,X_k)$. Then it follows from
    the fact that $X$ has no fixed time discontinuity that $X^N_k(W^N_k)$
    converges to $X_k(W_k)$. Using the assumption (A\ref{as:Feller}),
    this entails that $\widetilde{X}^N_{k+1}$ converges to a limiting
    process $\widetilde{X}_{k+1}$, which is distributed as $X$ started
    from $X_k(W_k)$. 

    Recalling that, by definition of the discrete spine under
    $\QQ^{k,N}_x$,
    \[
        X^N_{k+1} = \cat{X^N_k; W^N_k; \widetilde{X}^N_{k+1}},
    \]
    the claim is a consequence of the a.s.\ continuity of the
    concatenation map, which is proved in Lemma~\ref{lem:catContinuity}.
\end{proof}

\section{The recombination spine}
\label{S:1spine}

We now focus on the branching process with recombination. In this first
section we derive the properties of its $1$-spine. 

Using the formalism of the previous section, the branching process with
recombination can be constructed as a random marked tree, where the mark
space is the set of intervals of $\R$. 
According to the description of the branching process with recombination,
an individual with mark $I = [a,b]$ gives birth to $K(I)$ children, with
\[
    K(I) \sim \mathrm{Poisson}\big( 1 + \tfrac{\abs{I}}{N} \big).
\]
Then, each newborn experiences a recombination event with probability
\[
    r_N(I) \coloneqq \frac{2\abs{I}/N}{1+\abs{I}/N}.
\]
In the case of a recombination, the offspring inherits the interval $[a,
U]$ or $[U, b]$ with equal probability, where $U$ is uniformly
distributed over $I$. As in the previous section, we denote by
$\Xi(I)$ the offspring point process of a mother with interval $I$.
The objective of this section is to compute and characterize the
distribution $\QQ^{1,N}_I$ of the intervals along the $1$-spine 
and its large $N$ limit $\QQ^1_I$.

\paragraph{The $h$-transformed mark process.} 
Defining $\QQ^{1,N}_I$ first requires one to find an adequate harmonic
function for the branching process. In the branching process with
recombination, a simple calculation shows that the length of the
intervals is harmonic.

\begin{lemma}
    The function $h \colon I \mapsto \abs{I}$ is harmonic for the
    family of point processes $(\Xi(I))$. 
\end{lemma}

Let us now compute the distribution $\QQ^{1,N}_I$ of the $h$-transformed
process. According to \eqref{eq:harmonicKernel}, under $\QQ^{1,N}_I$, the
probability of experiencing no recombination in one time-step when
carrying interval $I$ is 
\[
    \Big( 1 + \frac{\abs{I}}{N} \Big) \big( 1-r_N(I) \big) = 1 -
    \frac{\abs{I}}{N}.
\]
When experiencing a recombination event, according to
\eqref{eq:harmonicKernel} the resulting interval is biased by $h$, that
is, biased by its length. This leads to the following description of the 
distribution of the intervals along the spine.

\begin{definition} \label{def:discreteSpine}
    The distribution $\QQ^{1,N}_I$ of the intervals along the $1$-spine 
    in the branching process with recombination is that of the
    discrete-time Markov chain $(I(n);\, n \ge 0)$ verifying $I(0) = I$,
    and conditional on $I(n) = [a,b]$, 
    \[
        I(n+1) = 
        \begin{dcases}
            [a,b] &\text{with probability $1 - \frac{b-a}{N}$}\\ 
            [a,a+U^*] &\text{with probability $\frac{b-a}{2N}$}\\ 
            [b-U^*,b] &\text{with probability $\frac{b-a}{2N}$}
        \end{dcases}
    \]
    where $U^*$ has the size-biased uniform distribution on 
    $[0, b-a]$. For convenience, $\QQ^{1,N}_R$ refers to 
    $\QQ^{1,N}_{[0,R]}$.
\end{definition}

\paragraph{Large $N$ convergence of the spine.}
As in the previous section, let $I^N$ denote the rescaled process
\[
    \forall t \ge 0,\quad I^N(t) = I\big( \floor{Nt} \big).
\]
We show that its large $N$ limit is given by the following process.

\begin{definition} \label{def:oneSpine}
    Let $\QQ^1_I$ denote the distribution of the continuous-time Markov
    process $(I(t);\, t \ge 0)$ started from $I$ and such that 
    jumps:
    \begin{itemize}
        \item from $[a,b]$ to $[a, a+U^*]$ at rate $(b-a)/2$;
        \item from $[a,b]$ to $[b-U^*, b]$ at rate $(b-a)/2$.
    \end{itemize}
    Again, $\QQ^1_R$ corresponds to $\QQ^1_{[0,R]}$.
\end{definition}

\begin{proposition} \label{prop:NconvergenceSpine}
    The process $(I^N(t);\, t \ge 0)$ under $\QQ^{1,N}_I$ converges
    in distribution for the Skorohod topology to $(I(t);\, t \ge 0)$
    under $\QQ^1_I$.
\end{proposition}

\begin{proof}
    The two processes $(I^N(t);\, t \ge 0)$ and $(I(t);\, t \ge 0)$ visit
    the same sequence of states, in distribution. Therefore, convergence
    in the Skorohod topology amounts to convergence of the jump times.
    Started from $[a,b]$, the time before the first jump of $(I^N(t);\, t
    \ge 0)$ is distributed as $T^N / N$, where $T^N$ is geometrically
    distributed with success probability $(b-a) / N$. It is clear that 
    $T^N / N$ converges in distribution to an exponentially distributed
    variable with mean $1/(b-a)$. Applying this convergence to the
    successive jump times of $(I^N(t);\, t \ge 0)$ readily proves the
    result.
\end{proof}

\subsection{Poisson construction of \texorpdfstring{$\QQ^1_R$}{Q\_R\^\ 1}}
\label{SS:poisson}

We are interested in the large $R$ properties of the spine. In this
section, we prove that the spine has a unique entrance law at infinity,
which can be constructed from a homogeneous Poisson point process. This
construction will also provide a coupling for the distribution of
$(I(t);\, t \ge 0)$ started from any initial condition.

First, for an interval $I = [a,b]$, it will be convenient to use the
notation
\[
    \lambda + \mu I = [\lambda + \mu a, \lambda + \mu b]
\]
for any reals $\lambda, \mu$.

Consider a homogeneous Poisson point process $P$ on $[0, \infty) \times \R$.
For any $t \ge 0$, consider the point process $P_t$ on $\R$ of atoms of $P$
with time coordinate in $[0, t]$ defined as 
\[
    \forall A,\quad P_t(A) = P([0,t] \times A).
\]
The atoms of $P_t(A)$ split the real line into infinitely many
subintervals. We are interested in the subinterval covering the origin.
More precisely, let $(x_i;\, i \in \Z)$ be the atoms of $P_t$,
labeled in such a way that
\[
    \ldots < x_{-1} < x_0 < 0 < x_1 < \ldots
\]
and define $I_P(t) = [x_0, x_1]$.

\begin{figure}
    \centering
    \includegraphics[width=.8\textwidth]{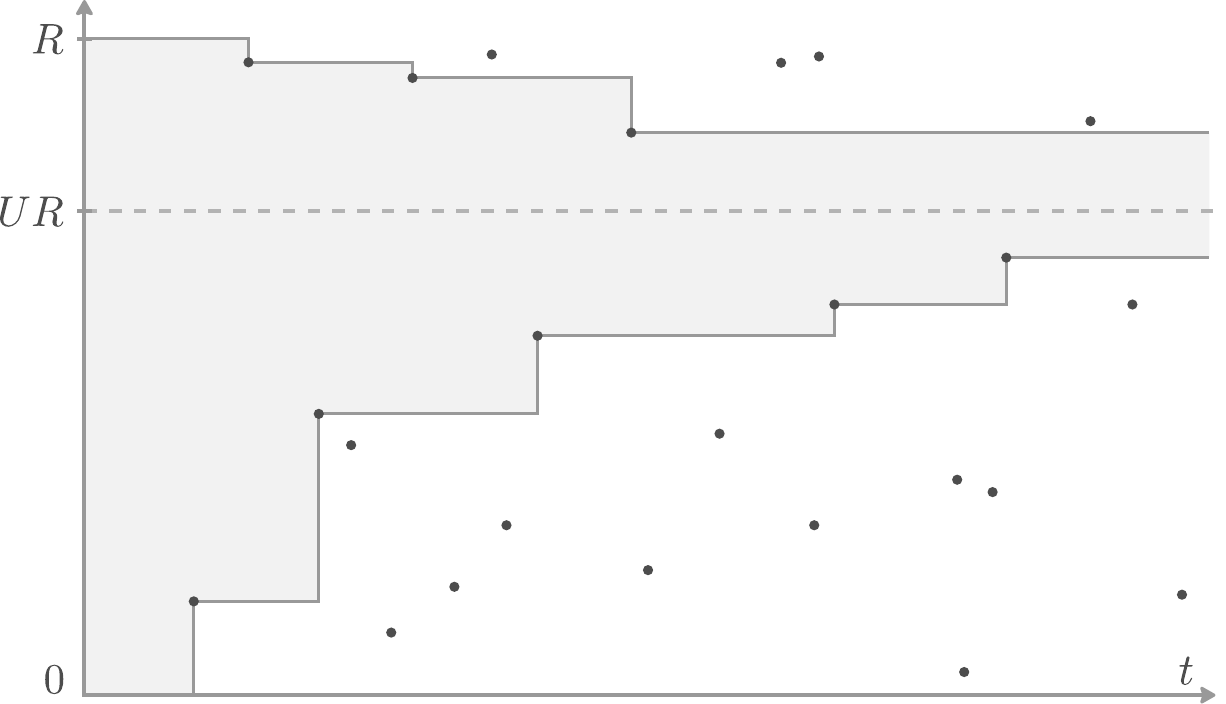}
    \caption{Illustration of the Poisson construction of $\QQ^1_R$.
    Atoms of $P$ are represented with dark circles. At each time $t$, 
    the vertical slice of the shaded region gives $I_R(t)$.}
\end{figure}

The following proposition shows that $(I_P(t);\, t \ge 0)$ corresponds to
the distribution of $(I(t);\, t \ge 0)$, started from infinity.

\begin{proposition} \label{prop:poissonConstruction}
    Let $M$ be uniformly distributed on $[0, 1]$ and independent of $P$.
    Then, for any $R \ge 0$, the process $(I_R(t);\, t \ge 0)$ defined as 
  
    \[
        \forall t \ge 0,\quad I_R(t) = MR  + I_P(t) \cap [-MR,(1-M)R]
    \]
    has distribution $\QQ^1_R$. Moreover, for any $t$, $M$ is uniformly
    distributed on $I_R(t)$.
\end{proposition}

The proposition will follow from the next simple result.

\begin{lemma} \label{lem:sampling}
    Let $U$ and $V$ be independent uniform r.v.\ on $[0, 1]$. Define 
    the interval 
    \[
        I = 
        \begin{cases}
            [0,U] &\text{if $V \le U$},\\
            [U,1] &\text{if $V > U$}
        \end{cases}
    \]
    Then $\abs{I}$ is independent of the event $\{V \le U\}$, $\abs{I}$
    is a size-biased uniform r.v.\ on $[0, 1]$, and $V$ is uniformly
    distributed on $I$.
\end{lemma}

\begin{proof}
    Let $A = \{ V \le U \}$ and $\bar{A} = \{ V > U \}$. For any test
    function $\phi$ and $\psi$, we can directly compute
    \begin{equation*}
        \P(\phi(V) \psi(\abs{I})\indic_{A} ) = \int_0^1 \int_0^u \phi(v)
        \psi(u) \diff v \diff u
        = \P(A) \int_0^1 2u \psi(u) \frac{1}{u} \int_0^u \phi(v) \diff v
        \diff u
    \end{equation*}
    and 
    \begin{align*}
        \P(\phi(V) \psi(\abs{I})\indic_{\bar{A}} ) &= \int_0^1 \int_u^1 \phi(v)
        \psi(1-u) \diff v \diff u \\
       &= \P(\bar{A}) \int_0^1 2(1-u) \psi(1-u) \frac{1}{1-u} \int_u^1 \phi(v) \diff v
        \diff u.
    \end{align*}
\end{proof}

\begin{proof}[Proof of Proposition~\ref{prop:poissonConstruction}]
    It is clear that $I_R(0) = [0, R]$.     
    Let us consider the sequence of jumps times $(T_i;\, i \ge 0)$ of
    $I_R$. By definition of the process, $T_{i+1}$ is the smallest time
    after $T_i$ such that there exists $X_{i+1}$ with $(T_{i+1}, X_{i+1})
    \in P$ and $X_{i+1} \in I_R(T_i)$. Then, if $I_R(T_i) = [a,b]$,
    \[
        I_R(T_{i+1}) = [a, X_{i+1}] \indic_{\{UR \le X_{i+1}\}}
    + [X_{i+1}, b] \indic_{\{X_{i+1} < UR\}}
    \]
    By the properties of homogeneous Poisson processes, $T_{i+1}-T_i$ is
    exponentially distributed with parameter $\abs{I_R(T_i)}$, and
    $X_{i+1}$ is uniformly distributed over $I_R(T_i)$.

    Using that $UR$ is uniformly distributed over $[0,R]$, a
    straightforward induction using Lemma~\ref{lem:sampling} proves that
    for any $i$, $UR$ is uniformly distributed on $I_R(T_i)$, and that
    $\abs{I_R(T_{i+1})}$ is a size-biased r.v.\ uniform on $I_R(T_i)$.
    This corresponds to the description of the transition mechanism of 
    $(I(t);\, t\ge 0)$ under $\QQ^1_R$ and proves the result.
\end{proof}

Define 
\[
    \forall t \ge 0,\quad X(t) = \abs{I(t)},\quad X_P(t) = \abs{I_P(t)}.
\]
The coupling provided by the previous representation can be used to
study the behavior of $(X(t);\, t \ge 0)$ as $R \to \infty$.

\begin{corollary} \label{cor:poissonLimit}
    As $R \to \infty$, the process $(X(t);\, t \ge 0)$ under $\QQ^1_R$
    converges in distribution to $(X_P(t);\, t \ge 0)$ for the
    topology of uniform convergence on every set $[\epsilon, \infty)$,
    $\epsilon > 0$. For any $t > 0$, $X_P(t)$ follows a
    $\mathrm{Gamma}(2, t)$ distribution, that is
    \[
        X_P(t) \sim t^2 x e^{-tx} \diff x.
    \]
\end{corollary}

\begin{proof}
    By the Poisson construction, if $X_R(t) \coloneqq \abs{I_R(t)}$, then
    $(X_R(t);\, t \ge 0)$ is distributed as $(X(t);\, t \ge 0)$ under
    $\QQ^1_R$. It is straightforward to see that for any $\epsilon > 0$, 
    for large enough $R$ we have $X_R(t) = X_P(t)$ for any $t \ge
    \epsilon$, proving the convergence part of the result.

    By well-known properties of Poisson point processes, $P_t$ is a
    homogeneous Poisson point process with rate $t$. Moreover both the
    first positive and negative atom of $P_t$ follow an exponential
    distribution with parameter $t$, and are independent. This proves
    that $X_P(t)$ is gamma distributed with the right parameters.
\end{proof}

\subsection{Self-similarity}
\label{SS:lamperti}
From the transition mechanism of $(I(t);\, t \ge 0)$, we see that $(X(t);\, t \ge
0)$ is also a Markov process. Under $\QQ^1_R$ it starts from $R$ and,
conditional on $X(t)$, it jumps at rate $X(t)$ to $X(t) U^*$, where $U^*$ is
a size-biased uniform r.v.\ on $[0, 1]$. As the jump rate of $(X(t);\, t \ge 0)$ 
at time $t$ is $X(t)$, the process  $(X(t);\, t \ge 0)$ is self-similar
with index $-1$, in the sense that for any constant $c > 0$, the
following identity in law holds
\[
    (c X(ct);\, t \ge 0) \text{ under $\QQ^1_R$ } \overset{(\mathrm{d})}{=} (X(t);\, t \ge 0)
    \text{ under $\QQ^1_{cR}$}.
\]

Positive Markov processes fulfilling the previous property are called
\emph{positive self-similar Markov processes} (pssMp), see
\cite{pardo13} for an introductory exposition.

\begin{remark}
    In this work we will not make use of this connection with pssMp since
    all computations can be carried out directly from the Poisson
    construction. However this link could be used to generalize our
    results to a larger class of branching processes on the intervals,
    with more general fragmentation rules for the offspring distribution.
\end{remark}

\subsection{Convergence of the rescaled spine}

Recall the definition of $F_R$ and $F_R^{-1}$ from \eqref{eq:rescaling}.
The following result provides the limit of the spine after rescaling time
according to $F_R^{-1}$. 
\begin{proposition}\label{prop:cvSpine}
    Let $0 < u_1 < \dots < u_n \le 1$. Then
    \[
        \big( F^{-1}_R(u_1)  X \circ F_R^{-1}(u_1), \dots, 
        F^{-1}_R(u_n)  X \circ F_R^{-1}(u_n) \big)
        \goesto[\QQ^{1}_R]{R \to \infty} (\gamma_1,\cdots, \gamma_n)
    \]
    where the $\gamma_i$'s are independent Gamma r.v.'s with parameter $(2,1)$.
\end{proposition}

\begin{proof}
    We show the result by induction on $n$. Recall that $X(F^{-1}_R(u))$ has
    the same distribution as $Y_1\wedge RU + Y_2 \wedge R(1-U)$ where
    $Y_1$, $Y_2$ are independent and exponentially distributed with mean
    $1/F^{-1}_R(u)$ and $U$ is uniform on $[0, 1]$. For every $u\in(0,1]$,
    $1/F^{-1}_R(u) \sim R^{1-u} = o(R)$ so that $Y_1\wedge RU + Y_2 \wedge
    R(1-U) = Y_1+Y_2$ with a probability going to $1$. This establishes
    the result at stage $1$.

    Let us now assume that the property is satisfied at stage $n$.
    Conditional on the process $X$ up to time  $F_R^{-1}(u_n)$, the Markov
    property implies that  the spine at $F^{-1}_R(u_{n+1})$ is distributed as
    $Y_1 \wedge R_n U + Y_2 \wedge R_n(1-U)$ where $R_n = X(F_R^{-1}(u_n))$,
    and $Y_1$, $Y_2$ are independent and exponentially distributed with mean
    $\big(F^{-1}_R(u_{n+1}) - F^{-1}_R(u_n)\big)^{-1}$. Since we have
    \[
        \frac{1}{F^{-1}_R(u_{n+1}) - F^{-1}_R(u_n)} 
        \sim \frac{1}{F^{-1}_R(u_{n+1})} 
        = o\Big( \frac{1}{F^{-1}_R(u_n)} \Big) = o( R_n ),
    \]
    as in the case $n = 1$, this implies that $F^{-1}_R(X(u_{n+1}))$ is
    converging to an independent $\gamma_{n+1}$ random variable.
\end{proof}

\section{The recombination \texorpdfstring{$k$}{k}-spine tree}
\label{S:kspine}

In the previous section we have characterized the large $N$, large $R$
behavior of the process giving the marks along a single spine. We now
provide a similar characterization for the $k$-spine tree. We start with
the following definition.

\begin{definition}\label{def:k-spine}
    Let us denote by $\QQ^k_I$ the law of some r.v.\ $(I_1,\dots,I_k)$
    and $(W_1,\dots, W_{k-1})$ such that:
    \begin{itemize}
        \item $(W_1,\dots, W_{k-1})$ are i.i.d.\ r.v.'s on $[0,1]$ with c.d.f. $F_R$;
        \item $I_1$ has distribution $\QQ^1_I$ and is independent of $(W_1,\dots, W_{k-1})$;
        \item $I_{i+1} = \cat{I_i; W_i; I'}$, where conditional on
            $(W_1,\dots, W_i)$ and $(I_1,\dots, I_i)$, $I'$ has
            distribution $\QQ^1_{I_i(W_i)}$.
    \end{itemize}
    We also use the shorter notation $\QQ^k_R$ for $\QQ^k_{[0,R]}$.
\end{definition}

\subsection{Convergence of the tree}

According to Proposition~\ref{prop:cvSpine}, it is natural to rescale
time using $F_R$ as follows.

\begin{definition}
    We consider the rescaled $k$-spine measure as
    \[
        \Qcal^k_R = \QQ^{k}_R \circ F_R^{-1},
    \]
    in the sense that under $\Qcal^k_R$
    \begin{itemize}
        \item The branch times are distributed as $(F_R(W_1), \dots, F_R(W_{k-1}))$;
        \item The spatial processes are distributed $(I_1 \circ F_R^{-1}, \dots, I_k \circ F_R^{-1})$;
    \end{itemize}
    where $(W_1,\dots, W_{k-1})$ and $(I_1,\dots, I_k)$ are distributed
    respectively as the branch times and the spatial processes under
    $\QQ^k_R$.
\end{definition}

The following result is a straightforward extension of
Proposition~\ref{prop:cvSpine}.

\begin{proposition}\label{prop:convTree}
    Under the rescaled $k$-spine measure $\Qcal_R^k$:
    \begin{enumerate}[(i)]
        \item The branch times $(W_1,\dots, W_{k-1})$ are distributed as
            i.i.d.\ uniform random variables on $[0,1]$.
        \item Conditional on the $W_i$'s
    \begin{multline*}
        \Big( \big( F_R^{-1}(W_1) X_1(W_1),\dots, 
                    F^{-1}_R(W_{k-1}) X_{k-1}(W_{k-1}) \big), 
              \big( X_1(1),\dots, X_{k}(1) \big) \Big) \\
        \goesto[\Qcal_R^{k}]{R\to\infty} 
        \big( (\gamma_1,\dots, \gamma_{k-1}), 
        (\bar{\gamma}_1,\dots, \bar{\gamma}_{k}) \big)
    \end{multline*}
    where the $\gamma_i$'s and $\bar \gamma_i$'s are independent Gamma
    r.v.'s with parameter $(1,2)$.
    \end{enumerate}
\end{proposition}

\begin{proof}
    The proof goes along the same line the one of
    Proposition~\ref{prop:cvSpine} and is left to the interested
    reader.
\end{proof}

\subsection{Convergence of the chromosomic distance}

In this section, we prove that the genealogical distance and the
rescaled chromosomic distance coincide in the large $R$ limit under 
$\mathcal{Q}^k_R$. For $i < j$, set 
\[
    W_{i,j} = W_{j,i} = \min \{ W_i, \dots, W_{j-1} \}
\]
to be the time at which branches $i$ and $j$ split. Define
\[
    \forall i,j \le k,\quad  d(i,j) = 1 - W_{i,j},
\]
which is the genealogical distance between the leaves of the $k$-spine
tree.

Conditional on $(I_1, \dots, I_k)$ and $(W_1, \dots, W_{k-1})$, let $M_i$
be uniformly distributed on $I_i(1)$. We define 
\[
    \forall i, j \le k,\quad D(i,j) = \abs{M_i - M_j}
\]
which is the chromosomic distance between the leaves. For later purpose,
we also introduce the corresponding rescaled distance,
\begin{equation} \label{eq:distances}
    \forall i,j \le k,\quad 
    \bar{d}_R(i,j) = 1 - F_R(W_{i,j}),\quad
    \bar{D}_R(i,j) = \frac{\log D(i,j) \vee 2}{\log R}.
\end{equation}

The next result provides an interesting relation between the genealogy of
the branching process and the ``geography'' along the chromosome. Namely,
on a logarithmic scale, the distance between two segments on the
chromosome is directly related to the genealogy of the two segments.

\begin{lemma} \label{lem:distanceConvergence}
    We have 
    \[
        \forall i,j \le k,\quad 
        \frac{\log D(i,j)}{d(i,j) \log R} 
        \goesto[\mathcal{Q}_R^k]{R \to \infty} 1.
    \]
\end{lemma}

\begin{proof}
    Let us work under $\QQ^k_R$ and let $i < j$. By construction of the
    $k$-spine, conditional on $I_i(W_{i,j})$, $(I_i(t+W_{i,j});\, t \ge 0)$
    and $(I_j(t+W_{i,j});\, t \ge 0)$ are independent and distributed as 
    $\QQ^1_{I_i(W_{i,j})}$. We know from the Poisson construction that,
    conditional on $I_i(W_{i,j})$, $M_i$ and $M_j$ are independent
    uniform variables on that interval. Therefore, 
    \[
        \frac{\abs{M_i-M_j}}{X_i(W_{i,j})}
    \]
    is a $\mathrm{Beta}(1,2)$ r.v.    
    Write
    \[
        \log D_{i,j} = \log \abs{M_i - M_j}
        = \log \frac{\abs{M_i - M_j}}{X_i(W_{i,j})}  
        + \log W_{i,j} X_i(W_{i,j})
        - \log W_{i,j}.
    \]
    From the previous point and Proposition~\ref{prop:cvSpine}, 
    \[
        \frac{1}{\log R}
        \log \frac{\abs{M_i - M_j}}{X_i(W_{i,j})}
        + \frac{\log W_{i,j} X_i(W_{i,j})}{\log R}
        \goesto[\QQ^k_R]{R \to \infty} 0.
    \]
    Moreover
    \[
        \frac{\log W_{i,j}}{\log R} = (F_R(W_{i,j})-1)(1+o_R(1))
    \]
    so that 
    \[
        \frac{\log D_{i,j}}{(1-F_R(W_{i,j})) \log R} 
        \goesto[\QQ^k_R]{R \to \infty} 1.
    \]
    The result follows by noting that the r.v.\ on the left-hand side
    has the same distribution under $\QQ^k_R$ as
    \[
        \frac{\log D_{i,j}}{d_{i,j} \log R} 
    \]
    under $\Qcal^k_R$.
\end{proof}

\subsection{Proof of the main result}
\label{SS:mainProof}

We can now proceed to the proof of our main result.

\begin{proof}[Proof of Theorem~\ref{thm:main2}]
    In order to ease the exposition, we only prove the result for $t =
    1$, but the proof is easily adapted for general $t > 0$.

    Recall that $\QQ^{1,N}_R$ denotes the distribution of the $1$-spine
    provided in Definition~\ref{def:oneSpine}. Let $\QQ^{k,N}_R$ be the
    corresponding $k$-spine distribution, with i.i.d.\ branch times 
    $(W_1, \dots, W_{k-1})$ such that 
    \begin{equation} \label{eq:branchTimes}
        W_i \overset{\mathrm{(d)}}{=} \floor{W N},\qquad
        \forall u \le 1,\quad \P(W \le u) = F_R(u).
    \end{equation}
    for the function $F_R$ defined in \eqref{eq:rescaling}.

    Set
    \[
        \forall i,j \le k, \quad 
        \bar{d}^N_R(i,j) = 1-F_R(W^N_{i,j}), \quad
        \bar{D}^N_R(i,j) = \frac{\log \abs{M_i - M_j} \vee 2}{\log R},
    \]
    where the $M_i$ are uniformly distributed on the $I_i(N)$.
    In order to use Theorem~\ref{thm:branchingConvergence}, we need to
    compute the limit of
    \[
            \QQ^{k,N}_R \Big[ \Delta_k \cdot 
                \phi\big( \bar{d}^N_R(i,j), \bar{D}^N_R(i,j), X^N_i(1);\,
                i, j \le k \big)
        \Big]
    \]
    where for the branching process with recombination,
    \[
        \Delta_k = \prod_{u \in S} \frac{1}{d_u!} \prod_{i=1}^{k-1}
        \frac{X_i^N(W^N_i)}{\delta_N F_R(W_i^N)} 
        \prod_{i=1}^k \frac{1}{X^N_i(1)}
    \]
    with
    \[
        \forall x \ge 0,\quad \delta_N F_R(x) = N \big({F_R(x+\tfrac{1}{N})- F_R(x)}\big).
    \]
    According to Proposition~\ref{prop:NconvergenceSpine}, under
    $\QQ^{1,N}_R$ the process $I^N$ converges to the limiting process
    with distribution $\QQ^1_R$ introduced in Definition~\ref{def:oneSpine}. 
    Therefore, Proposition~\ref{prop:Nconvergence} proves that 
    \[
        \big( (I^N_1,\dots,I^N_k), (W^N_1,\dots,W^N_{k-1}) \big) 
        \goesto[\QQ^{k,N}_R]{N \to \infty}
        \big( (I_1,\dots, I_k), (W_1,\dots, W_{k-1}) \big)
    \]
    where the limiting variables have distribution $\QQ^k_R$. Since the
    variables $(W_1,\dots,W_{k-1})$ are a.s.\ distinct under $\QQ^k_R$, it entails that 
    \[
        \prod_{u \in S} \frac{1}{d_u!} \goesto[\QQ^{k,N}_R]{N \to \infty}
        \frac{1}{2^{k-1}}.
    \]
    Corollary \ref{cor:negativeN} provides enough uniform integrability
    to conclude that
    \begin{align*}
        \lim_{N \to \infty} &\QQ^{k,N}_R \Big[ \prod_{u \in S} \frac{1}{d_u!} 
                \prod_{i=1}^{k-1}\frac{X_i^N(W^N_i)}{
                \delta_{N} F_R(W_{i}^N)} \prod_{i=1}^k
                \frac{1}{X^N_i(1)}  \cdot
                \phi\big( \bar{d}^N_R(i,j), \bar{D}^N_R(i,j), X^N_i(1);\, i,j \le k\big)\Big] \\
        &= \frac{1}{2^{k-1}} \QQ^k_R \Big[
                \prod_{i=1}^{k-1} \frac{X_i(W_i)}{F_R'(W_i)} \prod_{i=1}^k
                \frac{1}{X_i(1)}  \cdot
                \phi\big( \bar{d}_R(i,j), \bar{D}_R(i,j), X_i(1);\, i,j \le k\big)\Big] \\
        &= \frac{1}{2^{k-1}} {\cal Q}^k_R \Big[
                \prod_{i=1}^{k-1} X_i(W_i) (F^{-1}_R)'(W_i) 
                \prod_{i=1}^k \frac{1}{X_i(1)}  \cdot
                \phi\big( d(i,j), \bar{D}_R(i,j), X_i(1);\, i,j \le k\big)\Big]
    \end{align*}
    where the last line follows by rescaling time according to
    \eqref{eq:rescaling} and the distances are defined in
    \eqref{eq:distances}. The first observation is that, a.s., 
    \[
        \frac{1}{\log R} \frac{(F_R^{-1})'(W_i)}{F_R^{-1}(W_i)} =
        \frac{R^{W_i}}{R^{W_i}-1} \longrightarrow 1.
    \]
    The second observation is that conditional on $(W_1, \dots,
    W_{k-1})$, according to Proposition~\ref{prop:convTree}, 
    \[
        \big( (F^{-1}_R(W_i) X_i(W_i);\, i \le k-1), ( X_i(1);\, i \le k ) \big)
        \goesto[\Qcal^k_R]{R \to \infty}
        \big( (\gamma_1, \dots, \gamma_{k-1}), (\bar{\gamma}_1, \dots, \bar{\gamma}_k) \big)
    \]
    where the limiting r.v.\ are $\mathrm{Gamma}(2,1)$ distributed and
    independent. Let $\beta\in(1,2)$ and define
    \[
        \forall u \in [0,1], \quad \Phi(u) = \frac{1}{\log R \, F_R'(u) } = u +\frac{1}{R-1}
    \]
    By reverting the previous change of variable 
    \begin{equation*}
        \Qcal^k_R \Big[ 
            \prod_{i=1}^{k-1} \Big( \frac{(F^{-1}_R)'(W_i)X_i(W_i)}{\log R} \Big)^\beta 
            \prod_{i=1}^k \frac{1}{X_i(1)^\beta} \Big] 
         = \QQ^k_R \Big[
                \prod_{i=1}^{k-1} \big(X_i(W_i) \Phi(W_i)\big)^\beta \prod_{i=1}^k
           \frac{1}{X_i(1)^\beta} \Big].
   \end{equation*}
    By Corollary~\ref{cor:negativeN}, the r.h.s.\ is uniformly bounded in
    $R$. This provides enough uniform integrability to get
    \begin{align*}
    \lim_{R\to\infty} \Big(\frac{1}{2\log R}\Big)^{k-1} 
    \Qcal^k_R \Big[
            &\prod_{i=1}^{k-1} X_i(W_i) (F_R^{-1})'(W_i)
            \prod_{i=1}^k \frac{1}{X_i(1)}  \cdot
            \phi\big( d(i,j), \bar{D}_R(i,j), X_i(1);\, i,j \le k \big)
        \Big] \\
            &= \E\Big[ \prod_{i=1}^k \frac{1}{\gamma_i} \cdot 
                       \phi\big( 1-W_{i,j}, 1-W_{i,j}, \gamma_i;\, i,j \le k \big) \Big] \\
            &= \E\Big[ \phi\big( 1-W_{i,j}, 1-W_{i,j}, Y_i;\, i,j \le k \big) \Big]
    \end{align*}
    where $\gamma_i$ are i.i.d.\ $\mathrm{Gamma}(2,1)$ r.v., $Y_i$ are
    i.i.d.\ standard exponential random variable, and we have used
    Lemma~\ref{lem:distanceConvergence} for the convergence of the
    distances.

    We can now apply Theorem~\ref{thm:branchingConvergence} with 
    \[
        \alpha_N = \frac{1}{N \log R},\qquad
        \beta_N(x) = x,\qquad
        \gamma_N(t) = 1 - F_R\big(1- \tfrac{t}{N} \big).
    \]
    Recall that $Z_N$ denotes the population size and that 
    from Proposition~\ref{prop:survival} 
    \[
        \lim_{N \to \infty} N\PP_R(Z_N > 0) 
    \]
    exists and that 
    \[
        \lim_{R \to \infty} \lim_{N \to \infty} \frac{N\log R}{R}
        \PP_R(Z_N > 0) = 1.
    \]
    Therefore,
    \[
        \lim_{N \to \infty} \frac{k!N^{k-1}R}{(N\log R)^k\PP_R(Z_N > 0)}
            \QQ^{k,N}_R \Big[ \Delta_k \cdot \phi\big( \bar{d}^N_R(i,j),
                    \bar{D}^N_R(i,j), X_i^N(1);\, i,j \le k\big) \Big]
    \]
    exists and we have
    \begin{align*}
        \lim_{R \to \infty} 
        &\lim_{N \to \infty} \frac{k!N^{k-1}R}{(N\log R)^k\PP_R(\abs{T_N} > 0)}
            \QQ^{k,N}_R \Big[ \Delta_k \cdot \phi\big( \bar{d}^N_R(i,j), 
                    \bar{D}^N_R(i,j), X_i^N(1);\, i,j \le k\big) \Big] \\
        &= \lim_{R \to \infty} \lim_{N \to \infty} \frac{k!}{(\log R)^{k-1}}
            \QQ^{k,N}_R \Big[ \Delta_k \cdot \phi\big( \bar{d}^N_R(i,j), 
                    \bar{D}^N_R(i,j), X_i^N(1);\, i,j \le k\big) \Big] \\
        &=  k!\, \E\Big[ \phi\big( 1-W_{i,j}, 1-W_{i,j}, Y_i;\, i,j \le k \big) \Big]
    \end{align*}
    Applying Theorem~\ref{thm:branchingConvergence} for the large $N$
    limit, then Proposition~\ref{lem:convDetermining} for the large $R$ limit,
    and finally Proposition~\ref{prop:productSpace} and
    Proposition~\ref{prop:CPPpolynomials} for the polynomials of the
    Brownian CPP with independent marks proves the result. Note that
    since the total mass of the Brownian CPP is exponentially
    distributed, it fulfills the moment condition \eqref{eq:momentCondition}. 
    For the large $N$ limit, the size of the branching process with
    recombination is stochastically dominated by a Galton--Watson process
    with $\mathrm{Poisson}(1 + R/N)$ offspring distribution. The size of
    this process, conditional on survival at time $N$ and rescaled by
    $N$, is well-known to converge to an exponential distribution, see
    for instance \cite[Theorem~2.1]{oconnell95}. This readily shows that
    \eqref{eq:momentUMS} is also fulfilled by the large $N$ limit, for
    each fixed $R$.
\end{proof}

\begin{remark}
    In the above proof we have made use of Theorem~\ref{thm:convUMS} to
    identify the large $N$ limit of the genealogy. It is possible to prove the result
    without relying on our extension of the Gromov-weak topology by
    constructing the branching process with recombination by
    superimposing on a Galton--Watson tree with $\mathrm{Poisson}(1+R/N)$
    offspring distribution a process along the branches describing the
    recombination events as in \cite{baird_distribution_2003}. The
    large $N$ limit could then be expressed by means of the superprocess
    limit associated to that branching model.
\end{remark}

\paragraph{Acknowledgements.} We greatly thank Amaury Lambert for
initiating this project, for very helpful discussions, and for his
interest in this work.

\addcontentsline{toc}{section}{\refname}
\bibliographystyle{plain}
\bibliography{recombination.bib}

\appendix

\section{Appendix}

\subsection{Continuity of the concatenation map}

For two càdlàg functions $x$ and $y$ and a positive real $r$, recall that 
\[
    \cat{x; r; y} \colon t \mapsto 
    \begin{cases}
        x(t) &\text{if $t < r$}\\
        y(t-r) &\text{if $t \ge r$}
    \end{cases}
\]
denotes the concatenation of $x$ and $y$ at time $r$.

\begin{lemma} \label{lem:catContinuity}
    Suppose that we have two sequences of càdlàg functions $(x_n)$ and
    $(y_n)$ that converge to $x$ and $y$ respectively in the Skorohod
    topology. Moreover, suppose that $(r_n)$ converges to $r > 0$, and
    that $x$ is continuous at $r$. Then 
    \[
        \cat{x_n; r_n; y_n} \longrightarrow \cat{x; r; y}.
    \]
\end{lemma}

\begin{proof}
    Recall that a sequence $(z_n)$ converges to some $z$ in the Skorohod
    topology \tiff for each $T \ge 0$, we can find a sequence of 
    continuous increasing maps $\lambda_n$ such that 
    \begin{equation} \label{eq:sk1}
        \lim_{n \to \infty} \sup_{t \le T} \abs{\lambda_n(t) - t} = 0
    \end{equation}
    and 
    \begin{equation} \label{eq:sk2}
        \lim_{n \to \infty} \sup_{t \le T} \abs{z_n \circ \lambda_n(t) -
        z(t)} = 0,
    \end{equation}
    see for instance \cite[Proposition~5.3]{ethier_1986}.

    Fix $T > r$, and let $(\lambda_n)$ and $(\mu_n)$ fulfill
    \eqref{eq:sk1} and $\eqref{eq:sk2}$ with $(z_n)$ replaced by
    $(x_n)$ and $(y_n)$ respectively. For $\epsilon > 0$, let $\delta$ be
    such that 
    \[
        \abs{t-r} < \delta \implies \abs{x(t) - x(r)} < \epsilon.
    \]
    Using that 
    \[
        \lambda_n(r-\delta) \to r-\delta,\quad 
        \lambda_n(r+\delta) \to r+\delta, \quad 
        r_n \to r
    \]
    we know that for $n$ large enough $\lambda_n(r-\delta) < r_n <
    \lambda_n(r+\delta)$.
    Define 
    \[
        \forall t \ge 0,\quad \lambda'_n(t) = 
        \begin{dcases}
            \lambda_n(t) &\text{if $t \le r-\delta$}\\
            \lambda_n(r-\delta) + \big(t-(r-\delta)\big)\frac{r_n-\lambda_n(r-\delta)}{\delta} &\text{if $r-\delta < t \le r$}\\
            \lambda_n(r+\delta) + \big(t-(r+\delta)\big)\frac{\lambda_n(r+\delta)-r_n}{\delta} &\text{if $r < t \le r+\delta$}\\
            \lambda_n(t) &\text{if $t \ge r+\delta$}
        \end{dcases}
    \]
    which is continuous, increasing, and is designed so that
    $\lambda'_n(r) = r_n$ and coincides with $\lambda_n$ outside of
    $(r-\delta, r+\delta)$. Now,
    \begin{align*}
        \sup_{(r-\delta, r+\delta)} \abs{x_n \circ \lambda'_n(t) - x(t)}
        &\le 
        \sup_{(r-\delta, r+\delta)} \abs{x_n \circ \lambda'_n(t) - x(r)}
        + \sup_{(r-\delta, r+\delta)} \abs{x(r) - x(t)}\\
        &\le \sup_{(r-\delta, r+\delta)} \abs{x_n \circ \lambda_n(t) - x(r)}
        + \epsilon \\
        &\le \sup_{(r-\delta, r+\delta)} \abs{x_n \circ \lambda_n(t) - x(t)}
        + 2 \epsilon.
    \end{align*}
    Therefore, up to modifying $\lambda_n$ in this way and extracting a
    subsequence, we can always assume that, in addition to \eqref{eq:sk1} and 
    \eqref{eq:sk2}, we have $\lambda_n(r) = r_n$.

    Let us now denote 
    \[
        z_n = \cat{x_n; r_n; y_n},\quad z = \cat{x; r; y}
    \]
    and define
    \[
        \forall t \ge 0,\quad \nu_n(t) = 
        \begin{cases}
            \lambda_n(t) &\text{if $t < r$}\\
            \mu_n(t-r)+r_n &\text{if $t \ge r$}
        \end{cases}
    \]
    to be the concatenation of the two increasing maps. Clearly $(\nu_n)$
    fulfills \eqref{eq:sk1}, and using that $\nu_n(t) < r_n$ \tiff 
    $t < r$, we obtain 
    \[
        z_n \circ \nu_n(t) =
        \begin{cases}
            x_n \circ \lambda_n(t) &\text{if $t < r$}\\
            y_n \circ \mu_n(t-r) &\text{if $t \ge r$}
        \end{cases}
    \]
    from which it follows that
    \[
        \sup_{t \le T} \abs{z_n \circ \nu_n(t) - z(t)} \le 
        \sup_{t \le T} \abs{x_n \circ \lambda_n(t) - x(t)} 
        + \sup_{t \le T} \abs{y_n \circ \mu_n(t) - y(t)} 
    \]
    proving the result.
\end{proof}

\subsection{Uniform integrability results}

The aim of this technical section is to prove the uniform integrability
of the inverse of the $k$-spine.

\begin{lemma} \label{lem:negativeEstimate}
    For any $\alpha < 1$ and $t > 0$, there exists a constant $C$ such
    that 
    \[
        \forall R > 0,\, s \le t,\quad 
        \QQ^1_R \Big[ \frac{1}{X(s)^{1+\alpha}} \Big] \le
        C \Big( \frac{1}{R^{1+\alpha}} + 1 \Big).
    \]
    There also exists a constant $C'$ such that for $N$ sufficiently
    large
    \[
        \forall R > 0,\, s \le t,\quad 
        \QQ^{1,N}_R \Big[ \frac{1}{X^N(s)^{1+\alpha}} \Big] \le
        C' \Big( \frac{1}{R^{1+\alpha}} + 1 \Big).
    \]
\end{lemma}

\begin{proof}
    As $X^N$ and $X$ are non-increasing, we have 
    \[
        \QQ^1_R \Big[ \frac{1}{X(s)^{1+\alpha}} \Big] \le
        \QQ^1_R \Big[ \frac{1}{X(t)^{1+\alpha}} \Big]
    \]
    and 
    \[
        \QQ^{1,N}_R \Big[ \frac{1}{X^N(s)^{1+\alpha}} \Big] \le
        \QQ^{1,N}_R \Big[ \frac{1}{X^N(t)^{1+\alpha}} \Big].
    \]
    It is a direct consequence of the Poisson construction that 
    $X(t)$ under $\QQ^1_R$ is stochastically dominated by $X(t)$ under
    $\QQ^1_{R'}$ if $R \le R'$. Therefore 
    \[
        \forall R \ge 1,\quad 
        \QQ^1_R\Big[ \frac{1}{X(t)^{1+\alpha}} \Big] 
        \le
        \QQ^1_1\Big[ \frac{1}{X(t)^{1+\alpha}} \Big].
    \]
    Moreover, by self-similarity, 
    \[
        \forall R \le 1,\quad 
        \QQ^1_R\Big[ \frac{1}{X(t)^{1+\alpha}} \Big]
        = \frac{1}{R^{1+\alpha}} \QQ^1_1\Big[ \frac{1}{X(Rt)^{1+\alpha}}\Big]
        \le \frac{1}{R^{1+\alpha}} \QQ^1_1\Big[ \frac{1}{X(t)^{1+\alpha}} \Big].
    \]
    Together, these two identities yields that 
    \[
        \forall R > 0,\quad 
        \QQ^1_R\Big[ \frac{1}{X(t)^{1+\alpha}} \Big] 
        \le 
        \QQ^1_1\Big[ \frac{1}{X(t)^{1+\alpha}} \Big]
        \Big( \frac{1}{R^{1+\alpha}} + 1\Big).
    \]
    The first part of the result is proved provided that 
    \[
        \QQ^1_1\Big[ \frac{1}{X(t)^{1+\alpha}} \Big]  < \infty.
    \]
    Recalling that $X(t)$ has the same distribution as
    $Y_1 \wedge U + Y_2 \wedge (1-U)$ where $Y_1$ and $Y_2$ are
    exponentially distributed with mean $1/t$ and $U$ uniformly
    distributed on $[0,1]$, a direct computation shows that
    \[
        \QQ^1_1\Big[ \frac{1}{X(t)^{1+\alpha}} \Big]  
        \le \E\Big[ \Big(\frac{1}{Y_1+Y_2}\Big)^{1+\alpha} \Big]
        + 2 \E\Big[ \Big(\frac{1}{Y_1+U}\Big)^{1+\alpha} \Big]
        + 1 < \infty.
    \]

    The second part of the result follows from the fact that, for any 
    $c > 1$, for large enough $N$, 
    \[
        \forall R > 0,\quad X^N(t) \text{ under $\QQ^{1,N}_R$}
        \quad\overset{\mathrm{(d)}}{\ge}\quad
        X(ct) \text{ under $\QQ^1_R$}.
    \]
    To see this note that $X^N$ and $X$ visit the same sequence of states
    in distribution. Thus stochastic domination follows if $(X(ct);\, t
    \ge 0)$ jumps faster than $(X^N(t);\, t \ge 0)$. Started from $\rho$,
    the process $X^N$ jumps after a time $T^N / N$ whereas $X$ jumps
    after a time $T$ with
    \[
        T^N \sim \mathrm{Geometric}\Big( \frac{\rho}{N} \Big), \qquad
        T_\rho \sim \mathrm{Exponential}(\rho).
    \]
    A direct computation shows that, provided $1-\exp(-c\rho / N) \le \rho / N$,
    $T / c$ is stochastically dominated by $T^N / N$. The latter
    condition clearly holds for all $\rho \le R$ for large enough $N$.
\end{proof}

We now apply the previous estimate inductively to the $k$-spine, first to
be able to take the large $N$ limit, then to take the large $R$ limit.
Recall the notation $X^N_i$ for the rescaled marks along branch $i$ and
the notation $W^N_i$ for the rescaled $i$-th branch time, and the notation
\[
    \forall u \in [0,1],\quad \Phi(u) = u+ \frac{1}{R-1}.
\]

\begin{corollary} \label{cor:negativeN}
    For any $R > 1,k \ge 1$, and $\beta\in(1,2)$,
    \[
        \sup_{R > 1} \QQ^{k}_R\Big[ \prod_{i=1}^{k-1} (\Phi(W_i) X_{i}(W_i))^{\beta}
        \prod_{i=1}^k \Big(\frac{1}{X_i(1)}\Big)^{\beta}  \Big] 
        < \infty.
    \]
    Further,
    \[
        \sup_{N > 1} \QQ^{k,N}_R\Big[ 
            \prod_{i=1}^{k-1} \Big(\frac{X_{i}(W^N_i)}{\delta_N F_R(W_i^N)} \Big)^{\beta}
        \prod_{i=1}^k \frac{1}{X^N_i(1)^\beta}        
        \Big] < \infty.
    \]
\end{corollary}

\begin{proof}
We only prove the result in the continuum setting. The discrete case 
can be proved along the same lines.
    Applying the Markov property to $X_{k+1}$ at time $W_k$ and
    using Lemma~\ref{lem:negativeEstimate} yields
    \begin{align*}
        \QQ^{k+1}_R\Big[ 
        &\prod_{i=1}^{k} \big(\Phi(W_i) X_i(W_i)\big)^\beta \prod_{i=1}^{k+1} \Big(\frac{1}{X_i(1)}\Big)^{\beta} \Big]  \\
       &= \QQ^{k+1}_R\Big[ \prod_{i=1}^k 
            \big(\Phi(W_i)X_i(W_i)\big)^{\beta} \prod_{i=1}^k
            \frac{1}{X_i(1)^{\beta}}
            \QQ_{X_k(W_k)}^{1}\Big[
                \frac{1}{X(1-{W_k})^{\beta}}
            \Big] \Big]  \\
       &\le \QQ^{k+1}_R\Big[ \prod_{i=1}^k
            \big(\Phi(W_i)X_i(W_i)\big)^{\beta}  \prod_{i=1}^k
            \frac{1}{X_i(1)^{\beta}}
            C \Big(\frac{1}{X_k(W_k)^{\beta}} + 1 \Big)
         \Big] \\
       &\le \begin{multlined}[t]
      2^\beta C \QQ^{k}_R\Big[ \prod_{i=1}^{k-1}
            \big(\Phi(W_i)X_i(W_i)\big)^{\beta}
            \prod_{i=1}^k \frac{1}{X_i(1)^{\beta}}
         \Big] \\
       + C \QQ^{k}_R\Big[ 
           \big(\Phi(W_k)X_k(W_k)\big)^\beta \prod_{i=1}^{k-1}
            \big(\Phi(W_i)X_i(W_i)\big)^{\beta}
            \prod_{i=1}^k \frac{1}{X_i(1)^{\beta}}
       \Big]  \end{multlined}
    \end{align*}
    where in the last inequality, we used the fact that $\Phi(u)\leq 2$
    for every $u \in [0,1]$. Let $p, q \ge 0$ such that $\frac{1}{p} +
    \frac{1}{q} = 1$ and take $q$ close enough to $1$ such that $q \beta
    \in (1,2)$. By H\"older's inequality, the second term on the r.h.s.\ is
    bounded from above by
    \[
        \bigg( \QQ^{k}_R\Big[ 
            \prod_{i=1}^{k-1} \big(\Phi(W_i)X_i(W_i)\big)^{q\beta} 
            \prod_{i=1}^k \frac{1}{X_i(1)^{q\beta}}
        \Big] \bigg)^{1/q} 
        \bigg( \QQ^k_R\Big[
            \big(\Phi(W_1)X_1(W_1)\big)^{p\beta} 
        \big]\bigg)^{1/p}
    \]
    Further,
    \begin{align*}
      \QQ^k_R\Big[  \big(\Phi(W_1)X_1(W_1)\big)^{p\beta} \big] 
          &\le 2^{p\beta-1} \QQ^k_R\Big[ \big(W_1 X_1(W_1)\big)^{p\beta} \Big]
              + 2^{p\beta-1} \QQ^k_R\Big[ \Big(\frac{X_1(W_1)}{R-1}\Big)^{p\beta} \Big] \\
          &\le 2^{p\beta-1} \QQ^k_\infty\Big[ \big(W_1 X_1(W_1)\big)^{p\beta} \Big]
              + 2^{p\beta-1} \Big(\frac{R}{R-1}\Big)^{p\beta}
    \end{align*}
    Finally, $W_1 X_1(W_1)$ is a $\mathrm{Gamma}(2,1)$ r.v.\ under $\QQ^k_\infty$.
    The result follows by a straightforward induction since the case $k=1$
    was proved in Lemma~\ref{lem:negativeEstimate}.
\end{proof}

\end{document}